\documentclass[MSNbibl,number,citesort,seceqn,dvips]{arxbj}
\usepackage{mathbh}

\aid{0}
\volume{18}
\issue{2}
\pubyear{2012}
\firstpage{520}
\lastpage{551}
\doi{10.3150/11-BEJ353}

\makeatletter
\newcommand{\eqref}[1]{(\ref{#1})}

\newtheorem{theorem}{Theorem}[section]
\newtheorem{corollary}{Corollary}[section]
\newtheorem{lemma}{Lemma}[section]

\newtheorem{proposition}{Proposition}[section]
\newtheorem{prop}{Proposition}[section]

\newcommand{\R}{{\mathbb R}}
\newcommand{\Tr}{\operatorname{Tr}}
\newcommand{\eps}{{\varepsilon}}
\newcommand{\argmax}{\operatorname{argmax}}

\makeatother

\begin{document}
\begin{frontmatter}

\title{Central limit theorem and influence function for the MCD
estimators at general multivariate distributions}
\runtitle{CLT and influence function for the MCD}

\begin{aug}
%%%% inicialai - be tarpu
\author{\fnms{Eric A.} \snm{Cator}\thanksref{e1}\ead[label=e1,mark]{e.a.cator@tudelft.nl}}
\and
\author{\fnms{Hendrik P.} \snm{Lopuha\"a}\corref{}\thanksref{e2}\ead[label=e2,mark]{h.p.lopuhaa@tudelft.nl}}
\runauthor{E.A. Cator and H.P. Lopuha\"a}
\address{Delft University of Technology, DIAM, Mekelweg 4, 2628 CD
Delft, The Netherlands.\\
\printead{e1,e2}}
\end{aug}

% HISTORY:
\received{\smonth{6} \syear{2010}}
\revised{\smonth{12} \syear{2010}}

% ABSTRACT
%
\begin{abstract}
We define the minimum covariance determinant functionals for
multivariate location and scatter
through trimming functions and establish their existence at any
multivariate distribution.
We provide a precise characterization including a separating ellipsoid
property and prove
that the functionals are continuous.
Moreover, we establish asymptotic normality for both the location and
covariance estimator
and derive the influence function.
These results are obtained in a very general multivariate setting.
\end{abstract}

% KEYWORDS
%

\end{frontmatter}

%s1 #&#
\section{Introduction}\label{sec:intro}
Consider the minimum covariance determinant (MCD)
estimator introduced in \cite{rousseeuw85}, that is, for a sample
$X_1,X_2,\ldots,X_n$ from a distribution $P$ on $\R^k$ and
$0<\gamma\leq1$, consider subsamples $S\subset\{X_{1},\ldots,X_{n}\}$
that contain $h_n\geq\lceil{n\gamma}\rceil$ points. Define a
corresponding trimmed sample mean and sample covariance matrix by
%
%e1.1 #&#
\begin{eqnarray}\label{eq:def MCD estimator}
\widehat{T}_n(S) &=& \frac{1}{h_n}\sum_{X_i\in S}X_i,\nonumber\\ [-8pt]\\ [-8pt]
\widehat{C}_n(S) &=& \frac{1}{h_n}\sum_{X_i\in
S}\bigl(X_i-\widehat{T}_n(S)\bigr)\bigl(X_i-\widehat{T}_n(S)\bigr)'.\nonumber
\end{eqnarray}
Let $S_n$ be a subsample that minimizes $\det(\widehat{C}_n(S))$ over
all subsamples of size $h_n\geq\lceil{n\gamma}\rceil$, where
$\lceil
x\rceil$ denotes the smallest integer greater than or equal to
$x\in\R$. Then the pair $(\widehat{T}_n(S_n),\widehat{C}_n(S_n))$
is an
MCD estimator. Today, the MCD estimator is one of the most popular
robust methods to estimate multivariate location and scatter
parameters. These estimators, in particular the covariance estimator,
also serve as robust plug-ins in other multivariate statistical
techniques, such as principal component analysis
\cite{crouxhaesbroeck2000,serneelsverdonck2008}, multivariate linear
regression
\cite{agullocrouxvanaelst2008,rousseeuwvanaelstvandriessenagullo2004},
discriminant analysis \cite{hawkinsmclachlan1997}, factor analysis
\cite{pisonrousseeuwfilzmosercroux2003}, canonical correlations
\cite{taskinencrouxkankainen2006,zhou2009}, error-in-variables models
\cite{fekriruizgazen2004}, invariant coordinate selection
\cite{tylercritchleydumbgenoja2009}, among others (see also
\cite{hubertrousseeuwvanaelst} for a more extensive overview). For this
reason, the distributional and the robustness properties of the MCD
estimators are essential for conducting inference and perform robust
estimation in several statistical models.

The MCD estimators are known to have the same breakdown point as the
minimum volume ellipsoid estimators \cite{rousseeuw85}, and for a
suitable choice of $\gamma$ they possess the maximal breakdown point
possible for affine equivariant estimators (e.g., see
\cite{agullocrouxvanaelst2008,lopuhaarousseeuw1991}). However, their
asymptotic properties, such as the rate of convergence, limit
distribution and influence function, are not fully understood. Within
the framework of unimodal elliptically contoured densities, Butler,
Davies and Jhun \cite{butlerdaviesjuhn93} show that the MCD location
estimator converges at $\sqrt{n}$-rate towards a normal distribution
with mean equal to the MCD location functional. The rate of convergence
and limit distribution of the covariance estimator still remains an
open problem. Croux and Haesbroeck \cite{crouxhaesbroeck99} give the
expression for the influence function $\operatorname{IF}(x;C,P)$ of the MCD
covariance functional $C(P)$ at distributions $P$ with a unimodal
elliptically contoured density and use this to compute limiting
variances of the MCD covariance estimator. However, existence,
continuity and differentiability of the MCD functionals at perturbed
distributions is implicitly assumed, but not proven. Moreover, the
computation of the limiting variances via the influence function relies
on the von Mises expansion, that is,
%
%e1.2 #&#
\begin{equation}\label{eq:von mises cov}
\widehat{C}_n(S_n)-C(P)= \frac1n\sum_{i=1}^n
\operatorname{IF}(X_i;C,P) + \mathrm{o}_\mathbb{P}(n^{-1/2}),
\end{equation}
which has not been established. The distribution and robustness
properties of robust multivariate techniques that make use of the MCD,
depend on the distribution and robustness properties of the MCD
estimator, in particular those of the MCD covariance estimator. Despite
the incomplete asymptotic theory for the MCD, at several places in the
literature one prematurely assumes either a $\sqrt{n}$ rate of
convergence or asymptotic normality of the MCD covariance estimator, or
uses the influence function of the covariance MCD functional to
investigate the robustness of the specific multivariate method and to
determine limiting variances based on the heuristic \eqref{eq:von mises
cov}.

This paper is meant to settle these open problems and extend the
asymptotic theory for the MCD estimator in a very general setting that
allows a wide range of multivariate distributions. We will define the
MCD functional by means of trimming functions which are in a wide class
of measurable functions. Minimization of the determinant of our trimmed
covariance functional has similarities with minimizing the variation
functional corresponding to trimmed $k$-means considered in
\cite{gordaliza1991,cuesta-albertosgordalizamatran1995,cuesta-albertosgordalizamatranAS1997,garcia-escuderogordalizamatranAS1999}.
As opposed to our determinant functional, their variation functional
contains the trimming function in a linear way.
%This considerably facilitates compactness arguments that are used to
%establish existence and continuity for the functionals
%and which cannot be used in our setting.
However, once we have established the existence of our MCD functionals,
the characterization of our minimizing trimming functions is very
similar to the characterization found in
\cite{gordaliza1991,cuesta-albertosgordalizamatran1995}. In fact, part
of our characterization follows directly from results in
\cite{gordaliza1991}. These authors also recognized the advantage of
employing a flexible class of trimming functions, which allows a
uniform treatment at general probability measures, including empirical
measures and perturbed measures needed for our purposes. We believe
that obtaining our results for general multivariate distributions is an
important contribution of this paper. To justify this claim, we will
give several important examples of models where it is essential to
study the MCD estimator for a class of distributions that is wider than
the elliptically contoured distributions.

We prove existence of the MCD functional for any multivariate
distribution $P$ and provide a separating ellipsoid property for the
functional. Furthermore, we prove continuity of the functional, which
also yields strong consistency of the MCD estimators. Finally, we
derive an asymptotic expansion of the functional, from which we
rigorously derive the influence function, and establish a central limit
theorem for both MCD-estimators. We would like to emphasize that all
results are obtained under very mild conditions on~$P$ and that
essentially all conditions are satisfied for distributions with a
density. For distributions with an elliptically contoured density that
is unimodal we do not need any extra condition and recover the results
in \cite{butlerdaviesjuhn93} and \cite{crouxhaesbroeck99} as a special
case (see \cite{catorlopuhaa2009}).

The paper is organized as follows. In Section \ref{sec:def}, we define
the MCD functional for general underlying distributions, discuss some
of its basic properties and provide examples of models where it is
essential to study behavior of the MCD estimator for underlying
distributions that are beyond elliptically contoured distributions. In
Section \ref{sec:existence}, we prove existence of the MCD functional
and establish a separating ellipsoid property. Section
\ref{sec:continuity} deals with continuity of the MCD functionals and
consistency of the MCD estimators. Finally, in Section \ref{sec:asymp
norm} we obtain an asymptotic expansion of the MCD estimators and MCD
functional, from which we prove asymptotic normality and determine the
influence function. In order to keep things readable, all proofs and
technical lemmas have been postponed to an \hyperref[sec:appendix]{Appendix} at the end of the
paper.

%s2 #&#
\section{Definition}
\label{sec:def} Let $P$ be a probability measure on $\R^k$. To define
an MCD functional at $P$, we start by defining a trimmed mean and
trimmed covariance functional in the following way. For a~measurable
function $\phi:\R^k\to[0,1]$, define
%
%e2.1 #&#
\begin{eqnarray}\label{eq:def functionals}
T_P(\phi) &=&\frac1{\int\phi\,\mathrm{d}P}\int x \phi(x)
P(\mathrm{d}x),\nonumber\\[-8pt]\\ [-8pt]
C_P(\phi) &=& \frac1{\int\phi\,\mathrm{d}P}\int
\bigl(x-T_P(\phi)\bigr)\bigl(x-T_P(\phi)\bigr)'\phi(x) P(\mathrm{d}x).\nonumber
\end{eqnarray}
The function $\phi$ determines the trimming of the mean and covariance
matrix. For $\phi=1$, the above functionals are the ordinary mean and
covariance matrix corresponding to $P$. When $P=P_n$, the empirical
measure, and $\phi= \mathbh{1}_{S}$ for a subsample $S$, we recover
\eqref{eq:def MCD estimator}. Next, we fix a proportion $0<\gamma\leq
1$ and require $\phi$ to have at least mass $\gamma$, that is,
\[\label{eq:mass gamma}
\int\phi\,\mathrm{d}P\geq\gamma.
\]
To ensure that the functionals in \eqref{eq:def functionals} are well
defined, we take $\phi$ in the class
\[
K_P(\gamma) = \biggl\{\phi\dvtx\R^k\to[0,1]\dvt\phi\mbox{ measurable}, \int
\phi\,\mathrm{d}P
\geq\gamma, \int\|x\|^2\phi(x) P(\mathrm{d}x)<\infty\biggr\}.
\]
If there exists $\phi_P\in K_P(\gamma)$ which minimizes
$\det(C_P(\phi))$ over all $\phi\in K_P(\gamma)$, then the
corresponding pair
\[\label{eq:def MCD functional}
(T_P(\phi_P),C_P(\phi_P))
\]
is called an MCD functional at $P$. Note that, although for $\phi\in
K_P(\gamma)$ the functionals in~\eqref{eq:def functionals} are well
defined, the existence of a minimizing $\phi$ is not guaranteed.
Furthermore, if a minimizing $\phi$ exists, it need not be unique.

To complete our definitions, note that each trimming function $\phi$
determines an ellipsoid $E(T_P(\phi),C_P(\phi),r_P(\phi))$, where for
each $\mu\in\R^k$, $\Sigma$ symmetric positive definite, and $\rho>0$,
%
%e2.2 #&#
\begin{equation}\label{eq:def ellipsoid}
E(\mu,\Sigma,\rho) =  \{
x\in\R^k\dvt(x-\mu)'\Sigma^{-1}(x-\mu)\leq\rho^2  \},
\end{equation}
and
%
%e2.3 #&#
\begin{equation}\label{eq:def r_n}
r_P(\phi) = \inf\{ s>0 \dvt
P(E(T_P(\phi),C_P(\phi),s))\geq\gamma\}.
\end{equation}
If a minimizing trimming function $\phi_P$ exists, then
$E(T_P(\phi_P),C_P(\phi_P),r_P(\phi_P))$ is referred to as a
``minimizing'' ellipsoid.

Note that the functionals in \eqref{eq:def functionals} are affine
equivariant in the following sense. Fix a~non-singular $k\times k$
matrix $A$ and $b\in\R^k$ and let $h(x)=Ax+b$, for $x\in\R^k$. If
$X\sim P$, then $AX+b\sim Q=P\circ h^{-1}$. It is straightforward to
see that $\phi\in K_Q(\gamma)$ if and only if $\phi\circ h\in
K_P(\gamma)$, which yields
\[
T_Q(\phi)=A T_P(\phi\circ h)+b \quad\mbox{and}\quad C_Q(\phi)=A C_P(\phi
\circ
h)A',
\]
as well as $r_Q(\phi)=r_P(\phi\circ h)$. Furthermore, $\phi_Q$
minimizes $\det(C_Q(\phi))$ over $K_Q(\gamma)$ if and only if
$\phi_P=\phi_Q\circ h$ minimizes $\det(C_P(\phi))$ over
$K_P(\gamma)$. This means that if an MCD functional exists, it is
affine equivariant, that is, $T_Q(\phi_Q)=AT_P(\phi_P)+b$ and
$C_Q(\phi_Q)=AC_P(\phi_P)A'$.

Butler \textit{et al.} \cite{butlerdaviesjuhn93} define the MCD functional by
minimizing over all indicator functions~$\mathbh{1}_{B}$ of measurable
bounded Borel sets $B\subset\R^k$ with $P(B)=\gamma$. These indicator
functions form a subclass of $K_P(\gamma)$, that is sufficiently rich
when one considers unimodal elliptically contoured densities. However,
at perturbed distributions $P_{\eps,x}=(1-\eps)P+\eps\delta_x$, where
$\delta_x$ denotes the Dirac measure at $x\in\R^k$, their MCD
functional may not exist. Croux and Haesbroeck \cite{crouxhaesbroeck99}
solve this problem by minimizing over all functions
$\mathbh{1}_B+\delta\mathbh{1}_{\{x\}}$, with $x\notin B$ and
$P(B)+\delta P(\{x\})=\gamma$. These functions form a subclass of
$K_P(\gamma)$, that is sufficiently rich when one considers
single-point perturbations of unimodal elliptically contoured
densities, but the class $K_P(\gamma)$ allows for functions other than
$\mathbh{1}_B+\delta\mathbh{1}_{\{x\}}$ for which the determinant of
the covariance functional is strictly smaller. Moreover, minimization
over the more flexible class $K_P(\gamma)$ allows a uniform treatment
of the functionals in \eqref{eq:def functionals} at general\vadjust{\goodbreak} probability
measures, including measures with atoms. Important examples are the
empirical measure $P_n$ corresponding to a sample from $P$, in which
case the functionals relate to the MCD estimators, and perturbed
measures $P_{\eps,x}$, for which the functionals need to be
investigated in order to determine the influence function. It should be
noted that our Theorem \ref{th:characterization} does show that a
minimizer in the Croux--Haesbroeck sense does exist for all
distributions $P$, but this is not at all obvious before hand.

Definition \eqref{eq:def functionals} might suggest that minimization
of $\det(C_P(\phi))$ is hindered by the fact that the denominator
depends on $\phi$. However, the following property shows that if
a~minimum exists, it can always be achieved with a denominator in
\eqref{eq:def functionals} equal to $\gamma$. Its proof is
straightforward from definition \eqref{eq:def functionals}.
%
%l2.1 #&#
\begin{lemma}\label{lem:mass phi}
For any $0<\lambda\leq1$ and $\phi\in
K_P(\gamma)$, such that $\lambda\phi\in K_P(\gamma)$, we have
\[
T_P(\lambda\phi)=T_P(\phi),\qquad C_P(\lambda\phi)=C_P(\phi)
\quad\mbox{and}\quad r_P(\lambda\phi)=r_P(\phi).
\]
\end{lemma}

Since we can always construct a minimizing $\phi$ in such a way that
$\int\phi\,\mathrm{d}P=\gamma$, it is tempting to replace the term $\int\phi\,\mathrm{d}P$
in \eqref{eq:def functionals} by $\gamma$. However, we will not do
so, in order to keep enough flexibility for the functionals at
probability measures $P$ and trimming functions of the type
$\phi=\mathbh{1}_B$, for measurable $B\subset\R^k$ with
$P(B)>\gamma$.
An important example is the situation where $P$ is the empirical
measure.

%s2.1 #&#
\subsection{Examples of non-elliptical models where the MCD is relevant}
\label{subsec:examples} We will prove (see Theorem \ref{th:consistency
estimator}) that the MCD estimators converge (under mild conditions) to
the MCD functionals at $P$. These functionals might not be related in
any way to the expectation of $P$ or the covariance matrix (in fact,
our conditions allow for~$P$ whose expectation does not even exist) and
one might question the relevance of the MCD-functional for general $P$.

First of all, we believe that it is not unreasonable to consider the
MCD as a measure of location and scale on its own right, just like the
median and the MAD. Our results then show how the natural estimator of
this functional behaves. Especially in cases where the distribution has
a heavy tail, the MCD functional might provide more useful quantitative
information than the mean and covariance structure, for example for
confidence sets of future realizations of $P$. In~addition, we will
give some explicit examples in which it is very relevant to extend the
behavior of the MCD functional to general distributions.

%p2.1.0.1 #&#
\subsubsection*{Independent component analysis}
Consider a random vector $Z\in\R^k$ with a density $f$ that has the
property that for each coordinate, the mapping $y\mapsto
f(z_1,\ldots,y,\ldots,z_k)$ is a univariate, symmetric unimodal
function of $y$ for each fixed $z_1,\ldots, z_k$, and that $f$ is
invariant under coordinate-permutations. For example, this would be the
case if all the marginals of $f$ are independent and identically
distributed according to a univariate symmetric and unimodal
distribution. It is clear that if the MCD functional for $f$ is unique,
then from the symmetries it follows that the location functional is
zero, and the covariance functional is a constant times the identity
matrix. If we observe an affinely transformed sample from $f$, that is,
$X_1,\ldots,X_n$ where $X_i=AZ_i+\mu$ and $Z_i$ has density $f$, then
the MCD estimator would be a robust estimator of $\mu$ and $AA'$. Note
that the density of $X_1,\ldots,X_n$ is in general not elliptically
contoured. The uniqueness of the MCD functional for an $f$ of this kind
would be similar to the results in \cite{tatsuokatyler2006} for $S$-
and $M$-functionals. However, proving this is beyond the scope of this
paper, and might in fact be quite hard, given the depths of the results
in~\cite{tatsuokatyler2006}. The above example has close connections
with independent component analysis (ICA), a~highly popular method
within many applied areas, which routinely encounter multivariate data.
For a good overview see \cite{hyvarinenkarhunenoja2001}. The most
common ICA model considers~$X$ arising as a convolution of~$k$
independent components, that is, $X=AZ$, where~$A$ is nonsingular, and
the components of $Z$ are independent. The main objective of ICA is to
recover the mixing matrix $A$ so that one can `unmix' $X$ to obtain
independent components.

%p2.1.0.2 #&#
\subsubsection*{Invariant coordinate selection}
Invariant coordinate selection (ICS), recently proposed in
\cite{tylercritchleydumbgenoja2009}, compares two covariance estimators
through so-called ICS roots to reveal departures from an elliptically
contoured distribution. The authors suggest one of the covariance
estimators to be a class III scatter matrix, of which the MCD estimator
is an example. Determining whether ICS roots differ significantly, or
what power such a test would have, remains an open problem. This is
precisely where the distribution of the MCD estimator at elliptical
\textit{and} non-elliptical distributions is essential.

%p2.1.0.3 #&#
\subsubsection*{Contaminated distributions}
An important property for any robust estimator for location and scatter
is that it is able to recover to some extent the mean and covariance
matrix of the underlying distribution when this distribution is
contaminated. For instance, when the contamination has small total mass
or is very far away from the center of the underlying distribution, it
should not affect the corresponding functional too much. For our MCD
functional, this is precisely the content of the following theorem,
whose proof can be found in the \hyperref[sec:appendix]{Appendix}. These results rely heavily on
the methods used in this paper for general distributions, even if the
uncontaminated distribution $P$ is elliptically contoured.
%
%t2.1 #&#
\begin{theorem}\label{thm:contamination}
Let $P$ and $Q$ be two probability measures
on $\R^k$ and define for $x,r\in\R^k$ the translation $\tau_r(x)=x+r$.
Consider, for $\eps<1/2$, the mixture
\[
P_{r,\eps} = (1-\eps)P + \eps Q\circ\tau_r^{-1}.
\]
Denote by $\operatorname{MCD}_\gamma(\cdot)$ the MCD functional of level
$\gamma$. Choose $\gamma$ such that $\eps<\gamma<1-\eps$, and suppose
that
\[
P(H)<\frac{\gamma-\eps}{1-\eps}
\]
for all hyperplanes $H\subset\R^k$.
\begin{enumerate}[(ii)]
\item[(i)] Then
\[
\lim_{\eps\downarrow0} \operatorname{MCD}_\gamma(P_{r,\eps})
= \operatorname{MCD}_\gamma(P) \quad\mbox{and}\quad \lim_{\|r\|\to\infty}
\operatorname{MCD}_\gamma(P_{r,\eps}) = \operatorname{MCD}_{\gamma/(1-\eps)}(P),
\]
where the first limit should be interpreted as: every limit point is an
MCD functional at $P$ of level $\gamma$, and the second limit
similarly.
\item[(ii)] Furthermore, if in addition $Q$ has a bounded
support, then for all $\gamma\in(\eps,1-\eps)$, there exists
$r_0\geq
0$ such that
\[
\operatorname{MCD}_\gamma(P_{r,\eps}) = \operatorname{MCD}_{\gamma/(1-\eps)}(P)
\]
for all $r\in\R^k$ with $\|r\|\geq r_0$.
\end{enumerate}
\end{theorem}

As an illustration of Theorem \ref{thm:contamination}, consider an
elliptically contoured distribution $P$ with parameter $(\mu,\Sigma)$.
The second limit in (i) shows that if the contamination is far from
zero, the MCD functionals of the contaminated distribution are close to
$\mu$ and a~multiple of $\Sigma$. Part (ii) shows that for specific
types of contamination, for example, single point contaminations, the MCD
functionals at the contaminated distribution recovers these values
exactly. The proof of Theorem~\ref{thm:contamination} in principle
provides a constructive (but elaborate) way to find $r_0$ in terms of
$\eps, \gamma$, $P$ and the support of $Q$.

%s3 #&#
\section{Existence and characterization of an~MCD-functional}
\label{sec:existence} By definition, the matrix $C_P(\phi)$ is
symmetric non-negative definite. Without imposing any assumptions on
$P$, one cannot expect $C_P(\phi)$ to be positive definite. We will
assume that $P$ satisfies:
%
%e3.1 #&#
\begin{equation}\label{eq:hyperprop}
P(H)<\gamma\qquad \mbox{for every hyperplane $H\subset\R^k$.}
\end{equation}
This is a reasonable assumption, since if $P$ does not have this
property, then there exists a $\phi\in K_P(\gamma)$ with
$\det(C_P(\phi))=0$ (e.g., $\phi= 1_H$ with $P(H)\geq
\gamma$). This would prove the existence of a minimizing $\phi$, but
obviously the corresponding MCD-functional is not very useful.

We first establish the existence of a minimizing $\phi\in K_P(\gamma)$.
For later purposes, we do not only prove existence at $P$, but also at
probability measures $P_t$, for which the sequence $(P_t)$ converges
weakly to $P$, as $t\to\infty$. For ease of notation, we continue to
write $P_0$ instead of $P$ and for $t\geq0$ write
%
%e3.2 #&#
\begin{equation}\label{eq:def T C r}
T_t=T_{P_t},\qquad C_t=C_{P_t},\qquad r_t=r_{P_t}\quad
\mbox{and}\quad K_t(\gamma)=K_{P_t}(\gamma).
\end{equation}
The next proposition shows that eventually the smallest eigenvalue of
the covariance functional is bounded away from zero uniformly in $\phi$
and $t$.
%
%p1 #&#
\begin{prop}\label{prop:smalleig}
Suppose $P_0$ satisfies \eqref{eq:hyperprop} and
let $P_t\to P_0$ weakly. Then there exists $\lambda_0>0$ and $t_0\geq
1$ such that\vadjust{\goodbreak} for $t=0$, all $t\geq t_0$, all $\phi\in K_t(\gamma)$, and
all $a\in\mathcal{S}^k$ (the sphere in $\R^k$), we have
\[
\int\bigl(a'\bigl(x-T_t(\phi)\bigr)\bigr)^2\phi(x) P_t(\mathrm{d}x) \geq\lambda_0.
\]
In particular, this means that the smallest eigenvalue of $C_t(\phi)$
is at least $\lambda_0$.
\end{prop}

An immediate corollary is that if $\det(C_t(\phi))$ is uniformly
bounded, there exists a~compact set that contains the location and
covariance functionals for sufficiently large $t$ (see Lemma
\ref{lem:compact set} in the \hyperref[sec:appendix]{Appendix}). This will become very useful in
establishing continuity of the functionals in Section
\ref{sec:continuity}. For the moment, we use this result to show that
for minimizing $\det(C_t(\phi))$, one may restrict to functions $\phi$
with bounded support.

For $R>0$, define the ball $B_R = \{x\in\R^k \dvt \|x\|\leq R\}$ and for
$t\geq0$ define the class
\[
K_t^R(\gamma) = \bigl\{ \phi\in K_t(\gamma)\dvt \{\phi\ne0\}\subset B_R\bigr\}.
\]
Clearly, $K^R_t(\gamma) \subset K_t(\gamma)$. The next proposition
shows that for any $\phi\in K_t(\gamma)$ we can always find a $\psi$
with bounded support in $K_t^R(\gamma)$ that has a smaller determinant.
%
%p2 #&#
\begin{prop}\label{prop:boundedsupp}
Suppose that $P_0$ satisfies
\eqref{eq:hyperprop} and $P_t\to P_0$ weakly. There exists $R>0$ and
$t_0\geq1$ such that for $t=0$, all $t\geq t_0$ and all $\phi\in
K_t(\gamma)$, there exists $\psi\in K^R_t(\gamma)$ with
\[
\det(C_t(\psi))\leq\det(C_t(\phi)).
\]
\end{prop}

Proposition \ref{prop:boundedsupp} illustrates the general heuristic
that if $\phi$ has $P_t$-mass far away from~$T_t(\phi)$, then moving
this mass closer towards $T_t(\phi)$ will decrease the determinant of
the covariance matrix. Together with Proposition \ref{prop:smalleig}
this establishes the existence of at least one MCD functional for the
probability measure $P_0$. Moreover, if $P_t\to P_0$ weakly, then at
least one MCD functional exists for $P_t$ for sufficiently large $t$.

%t3.1 #&#
\begin{theorem}
\label{th:existence} Suppose $P_0$ satisfies \eqref{eq:hyperprop} and
let $P_t\to P_0$ weakly. Then there exists $R>0$ and $t_0\geq1$, such
that for $t=0$ and $t\geq t_0$, there exists $\phi_t\in K_t^R(\gamma)$,
which minimizes $\det(C_t(\phi))$ over $K_t(\gamma)$.
\end{theorem}

In the remainder of this section, we provide a characterization of a
minimizing $\phi$, which includes a separating ellipsoid property for
the MCD functional. A similar result has been obtained in
\cite{butlerdaviesjuhn93} for the empirical measure and in
\cite{crouxhaesbroeck99} for single-point perturbations of
distributions with a unimodal elliptically contoured density. Our
characterization is for general distributions and is very similar to
the characterizations obtained for trimmed means in
\cite{gordaliza1991} and \cite{cuesta-albertosgordalizamatran1995}. In
fact, the first part of our result follows directly from
\cite{gordaliza1991}. We will denote the interior of a set $E$ by
$E^\circ$, and the (topological) boundary by $\partial E$.

%t3.2 #&#
\begin{theorem}
\label{th:characterization} Let $\phi\in K_P(\gamma)$ be such that
$(T_P(\phi),C_P(\phi))$ is an MCD functional at $P$ and let
$E_P(\phi)=E(T_P(\phi),C_P(\phi),r_P(\phi))$ be the corresponding
minimizing ellipsoid. Then
\[
\int\phi\,\mathrm{d}P = \gamma\quad \mbox{and}\quad \mathbh{1}_{E_P(\phi)^\circ} \leq
\phi\leq\mathbh{1}_{E_P(\phi)},\qquad P\mbox{-a.e.}
\]
Furthermore, either $\phi=0$ on $\partial E_P(\phi)$ ($P$-a.e.), or
$\phi=1$ on $\partial E_P(\phi)$ ($P$-a.e.), or there exists $x\in
\partial E_P(\phi)$ such that $P(\partial E_P(\phi)) = P(\{x\})$.
\end{theorem}

The theorem shows that a minimizing trimming function $\phi$ is almost
the indicator function of an ellipsoid with center $T_P(\phi)$ and
covariance structure $C_P(\phi)$. When $P$ has no mass on the boundary
of the ellipsoid, then $\phi$ is equal to the indicator function of
this ellipsoid. If the interior of the ellipsoid
$E(T_P(\phi),C_P(\phi),r_P(\phi))$ has mass strictly smaller than
$\gamma$, then either $\phi$ equals~$1$ on the entire boundary of the
ellipsoid, in which case the (closed) ellipsoid has $P$ mass exactly~$\gamma$, or $P$ only has mass in exactly one point on the boundary,
and $\phi$ adapts its value in that point such that it has total $P$
mass $\gamma$.

Theorem \ref{th:characterization} holds for any probability measure, in
particular for the empirical measure~$P_n$ and for perturbed measures
$P_{\eps,x}$, in which case we obtain results analogous to Theorem 2 in
\cite{butlerdaviesjuhn93} and Proposition 1 in
\cite{crouxhaesbroeck99}, respectively. However, note that the
characterization in Theorem \ref{th:characterization} is more precise
and is such that the center and covariance structure of the separating
ellipsoid are exactly the MCD functionals themselves.

%s4 #&#
\section{Continuity of the MCD functional}
\label{sec:continuity} Consider a sequence $
(T_{t}(\phi_t),C_{t}(\phi_t) )$ of MCD functionals corresponding to a
sequence of probability measures $P_t\to P_0$ weakly. We investigate
under what conditions $ (T_{t}(\phi_t),C_{t}(\phi_t) )$ converges and
whether each limit point will be an MCD functional corresponding to
$P_0$. Our approach requires $\int\phi\,\mathrm{d}P_t\to\int\phi\,\mathrm{d}P_0$ uniformly
in minimizing $\phi$. The following condition on $P_0$ suffices:
%
%e4.1 #&#
\begin{equation}\label{eq:unif conv ellipsoid}
\sup_{E\in
\mathcal{E}}|P_t(E)-P_0(E)|\to0,\qquad \mbox{as }t\to\infty,
\end{equation}
where $\mathcal{E}$ denotes the class of all ellipsoids. This may seem
restrictive, but it is either automatically fulfilled for sequences
that are important for our purposes or a mild condition on $P_0$
suffices. For instance, when $P_t$ is a sequence of empirical measures,
then~\eqref{eq:unif conv ellipsoid} holds automatically by standard
results from empirical process theory (e.g., see Theorem~II.14 in
\cite{pollard84}) because the ellipsoids form a class with polynomial
discrimination or a Vapnik--Cervonenkis class. Condition \eqref{eq:unif
conv ellipsoid} also holds for sequences of perturbed measures
$P_{\eps,x}$, as $\eps\downarrow0$. In general, if $P_0(\partial C)=0$
for all measurable convex $C\subset\R^k$, then condition \eqref{eq:unif
conv ellipsoid} holds for any sequence $P_t\to P_0$ weakly (see Theorem
4.2 in \cite{rangarao62}). Note that this is always trivially true if
$P_0$ has a density.

For later purposes, we prove continuity not only for MCD functional
minimizing functions $\phi_t$, but for any sequence of functions
$\psi_t$ with uniformly bounded support that satisfy the same
characteristics as $\phi_t$ and for which $\det(C_t(\psi_t))$ is close
to $\det(C_t(\phi_t))$.
%
%t4.1 #&#
\begin{theorem}
\label{th:continuity} Suppose $P_0$ satisfies \eqref{eq:hyperprop}. Let
$P_t\to P_0$ weakly and suppose that \eqref{eq:unif conv ellipsoid}
holds. For $t\geq1$, let $\psi_t\in K_t(\gamma)$ such that $\psi
_t\leq
\mathbh{1}_{E_t}$, where $E_t=E(T_t(\psi_t),C_t(\psi_t),r_t(\psi_t))$,
and suppose there exist $R>0$ such that $\{\psi_t\ne0\}\subset B_R$,
for $t$ sufficiently large. Suppose that
\[
\det(C_t(\psi_t))-\det(C_t(\phi_t))\to0,\qquad \mbox{as
}t\to\infty,\vadjust{\goodbreak}
\]
where $\phi_t$ minimizes $\det(C_t(\phi))$ over $K_t(\gamma)$. Then
\begin{enumerate}[(ii)]
\item[(i)] there exist a convergent subsequence $
(T_{t_m}(\psi_{t_m}),C_{t_m}(\psi_{t_m}) )$;
\item[(ii)] the
limit point of any convergent subsequence is an MCD functional at
$P_0$.
\end{enumerate}
\end{theorem}

An immediate corollary is that in case the MCD functional at $P_0$ is
uniquely defined, all possible MCD functionals at $P_t$ are consistent.
For later purposes, we also need that~$r_t(\phi_t)$ converges. This may
not be the case if $P_0$ has no mass directly outside the boundary of
its minimizing ellipsoid. For this reason, we also require that
%
%e4.2 #&#
\begin{equation}\label{eq:cond boundary}
P_0\bigl(E\bigl(T_0(\phi_0),C_0(\phi_0),r_0(\phi_0)+\eps\bigr)\bigr)>\gamma\qquad \mbox{for all
}\eps>0.
\end{equation}
Note that this condition is trivially true if $P_0$ has a positive
density in a neighborhood of the boundary of
$E(T_0(\phi_0),C_0(\phi_0),r_0(\phi_0))$.

%c4.1 #&#
\begin{corollary}
\label{cor:consistency functional} Suppose $P_0$ satisfies
\eqref{eq:hyperprop} and that the MCD functional $
(T_0(\phi_0),C_0(\phi_0) )$ is uniquely defined at $P_0$. Let $P_t\to
P_0$ weakly and suppose that \eqref{eq:unif conv ellipsoid} holds. For
$t\geq1$, let $\psi_t\in K_t(\gamma)$ such that $\psi _t\leq
\mathbh{1}_{E_t}$, where $E_t=E(T_t(\psi_t),C_t(\psi_t),r_t(\psi_t))$,
and suppose there exist $R>0$ such that $\{\psi_t\ne0\}\subset B_R$,
for $t$ sufficiently large. Suppose that
\[
\det(C_t(\psi_t))-\det(C_t(\phi_t))\to0,
\]
where $\phi_t$ minimizes $\det(C_t(\phi))$ over $K_t(\gamma)$. Then,
\begin{enumerate}[(ii)]
\item[(i)] $ (T_{t}(\psi_t),C_{t}(\psi_t) )\to (T_0(\phi
_0),C_0(\phi_0) )$.
\end{enumerate}
If in addition $P_0$ satisfies \eqref{eq:cond boundary}, then
\begin{enumerate}[(ii)]
\item[(ii)] $r_t(\psi_t)\to r_0(\phi_0)$.
\end{enumerate}
\end{corollary}

Uniqueness of the MCD functional has been proven in
\cite{butlerdaviesjuhn93} for distributions $P_0$ that have a unimodal
elliptically contoured density. For general distributions, one cannot
expect such a general result. For instance, for certain bimodal
distributions or for a~spherically symmetric uniform distribution which
is positive on a large enough disc, the MCD functional is no longer
unique.

%s4.1 #&#
\subsection{Consistency of the MCD estimators}
\label{subsec:consistency MCD estimator} For $n=1,2,\ldots,$ let $P_n$
denote the empirical measure corresponding to a sample from~$P_0$. From
definitions \eqref{eq:def MCD estimator} and \eqref{eq:def
functionals}, it is easy to see that the MCD estimators can be written
in terms of the MCD functional as follows
%
%e4.3 #&#
\begin{eqnarray}\label{eq:rel est-func}
\widehat{T}_n(S_n) &=&
T_n(\mathbh{1}_{S_n}),\nonumber\\ [-8pt]\\ [-8pt]
\widehat{C}_n(S_n) &=& C_n(\mathbh{1}_{S_n}),\nonumber
\end{eqnarray}
where we use the notation introduced in \eqref{eq:def T C r}. Moreover,
define $\widehat{r}_n(S_n)=r_n(\mathbh{1}_{S_n})$. We should
emphasize\vadjust{\goodbreak}
that $\widehat{T}_n(S_n)$ and $\widehat{C}_n(S_n)$ may differ from the
actual MCD functionals $T_n(\phi_n)$ and $C_n(\phi_n)$, where $\phi_n$
minimizes the MCD functional for $P_n$. Obviously, if these differences
tend to zero, then consistency of the MCD estimators would follow
immediately from Theorem \ref{th:continuity}, but unfortunately we have
not been able to find an easy argument for this. However, we can show
that the determinants of the covariance matrices are close with
probability one, which suffices for our purposes.

%p4.1 #&#
\begin{proposition}
\label{prop:diff det subsample-phi} Suppose $P_0$ satisfies
\eqref{eq:hyperprop}. Then for each MCD estimator minimizing subsample
$S_n$ and each MCD functional minimizing function $\phi_n$, we have
\[
\det(\widehat{C}_n(S_n))-\det(C_n(\phi_n))=\mathrm{O}(n^{-1})
\]
with probability one.
\end{proposition}

This does not necessarily mean that $\widehat{T}_n(S_n)-T_n(\phi_n)$
and $\widehat{C}_n(S_n)-C_n(\phi_n)$ are also of the order $\mathrm{O}(n^{-1})$.
But in view of Corollary \ref{cor:consistency functional}, it suffices
to establish a separating ellipsoid property and uniform bounded
support for the minimizing subsample. The latter result can be found in
the \hyperref[sec:appendix]{Appendix}, whereas the separating ellipsoid property is stated in
the next proposition.

%p4.2 #&#
\begin{proposition}
\label{prop:characterization estimator} Let $S_n$ be a minimizing
subsample for the MCD estimator and define corresponding ellipsoid
$\widehat{E}_n=E(\widehat{T}_n(S_n),\widehat{C}_n(S_n),\widehat{r}_n(S_n))$.
Then $S_n$ has exactly $\lceil n\gamma\rceil$ points, $S_n\subset
\widehat{E}_n$ and $\widehat{E}_n$ only contains points of $S_n$.
\end{proposition}

This separating ellipsoid property is somewhat different from the one
in Theorem \ref{th:characterization} (for the empirical measure) and
from the one in \cite{butlerdaviesjuhn93}. The ellipsoid
$\widehat{E}_n$ has the MCD estimators as center and covariance
structure instead of the trimmed sample mean and covariance
corresponding to the minimizing subsample excluding a point that is
most outlying (see \cite{butlerdaviesjuhn93}). The advantage of the
characterization given in Proposition \ref{prop:characterization
estimator} is that integrating over $S_n$ or $\widehat{E}_n$ with
respect to $P_n$ is the same, which will become very useful later on.
We now have the following theorem.
%
%t4.2 #&#
\begin{theorem}
\label{th:consistency estimator} Suppose $P_0$ satisfies
\eqref{eq:hyperprop} and that the MCD functional $
(T_0(\phi_0),C_0(\phi_0) )$ is uniquely defined at $P_0$. For $n\geq1$,
let $S_n$ be a minimizing subsample for the MCD estimator. Then
\begin{enumerate}[(ii)]
\item[(i)] $ (\widehat{T}_n(S_n),\widehat{C}_n(S_n) )\to (T_0(\phi
_0),C_0(\phi_0) )$, with probability one.
\end{enumerate}
If, in addition $P_0$ satisfies \eqref{eq:cond boundary}, then
\begin{enumerate}[(ii)]
\item[(ii)] $\widehat{r}_n(S_n)\to r_0(\phi_0)$, with probability one.
\end{enumerate}
\end{theorem}

As a special case, where $P_0$ has a unimodal elliptically contoured
density, we recover Theorem 3 in \cite{butlerdaviesjuhn93}. With
Theorems \ref{th:continuity} and \ref{th:consistency estimator}, it
turns out that the difference between the MCD estimator $
(\widehat{T}_n(S_n),\widehat{C}_n(S_n) )$ and the MCD functional $
(T_{n}(\phi_n),C_{n}(\psi_n) )$ indeed tends to zero with probability
one. However, we were not able to find an easier, direct argument.

%BEJ353_Section5.tex
%s5 #&#
\section{Asymptotic normality and influence function}\vspace*{3pt}
\label{sec:asymp norm} For $n=1,2,\ldots,$ let $S_n$ be a minimizing
subsample for the MCD estimator and for ease of notation, write
\[
\label{eq:def MCD estimator brief}
\widehat\mu_n=\widehat{T}_n(S_n),\qquad
\widehat\Sigma_n=\widehat{C}_n(S_n)=\widehat{\Gamma}_n^2,\qquad
\widehat{\rho}_n=\widehat{r}_n(S_n)\quad \mbox{and}\quad
\widehat{E}_n=E(\widehat{\mu}_n,\widehat{\Sigma}_n,\widehat{\rho}_n),
\]
and define
$\widehat{\theta}_n=(\widehat{\mu}_n,\widehat{\Gamma}_n,\widehat
{\rho}_n)$
in $\R^k\times\operatorname{PDS}(k)\times\R$, where $\operatorname{PDS}(k)$
denotes the class of all positive definite symmetric matrices of order
$k$. Note that $\widehat{\Gamma}_n$ is uniquely defined in
$\operatorname{PDS}(k)$. Similarly, let $P_n$ denote the empirical measure
corresponding to a sample from $P_0$, and for $n=0,1,2,\ldots,$ let
$\phi_n$ be a minimizing trimming function for the MCD functional and
write
%
%e5.1 #&#
\begin{equation}\label{eq:def mu Gamma rho}
\mu_n=T_n(\phi_n),\quad
\Sigma_n=C_n(\phi_n)=\Gamma_n^2,\quad \rho_n=r_n(\phi_n)\quad
\mbox{and}\quad E_n=E(\mu_n,\Sigma_n,\rho_n),
\end{equation}
where $T_n$, $C_n$ and $r_n$ are defined in \eqref{eq:def T C r}, and
write $\theta_n=(\mu_n,\Gamma_n,\rho_n)$. According to Corollary
\ref{cor:consistency functional} and Theorem \ref{th:consistency
estimator}, under very mild conditions on $P_0$, we have
$\widehat{\theta}_n\to\theta_0$ and $\theta_n\to\theta_0$ with
probability one, where $\theta_0=(\mu_0,\Gamma_0,\rho_0)$ corresponds
to $P_0$ as defined in~\eqref{eq:def mu Gamma rho}. The limit
distribution of $\widehat{\theta}_n-\theta_0$ and $\theta_n-\theta_0$
are equal and can be obtained by the same argument. We briefly sketch
the main steps for the MCD estimator.

Consider the estimator matrix equation in \eqref{eq:def MCD estimator},
\[
\widehat{\Sigma}_n=\frac{1}{P_n(S_n)}\int_{S_n}(x-\widehat{\mu
}_n)(x-\widehat{\mu}_n)'
P_n(\mathrm{d}x).
\]
After multiplying from the left and the right by
$\widehat{\Gamma}_n^{-1}$, rearranging terms and replacing~$S_n$ by
$\widehat{E}_n$ (which leaves the integral unchanged according to
Proposition \ref{prop:characterization estimator}), we obtain
a~covariance valued $M$-estimator type score equation:
\[
0= \int_{\widehat{E}_n} \bigl(
\widehat{\Gamma}_n^{-1}(x-\widehat{\mu}_n)(x-\widehat{\mu
}_n)'\widehat{\Gamma}_n^{-1}-I_k
\bigr) P_n(\mathrm{d}x).
\]
Similarly, one can obtain a vector valued $M$-estimator type score
equation from the location equation in \eqref{eq:def MCD estimator} and
the equality $P_n(\widehat{E}_n)=\lceil
n\gamma\rceil/n=\gamma+\mathrm{O}(n^{-1})$ can be put into a real valued score
equation. Putting everything together, we conclude that
$\widehat{\theta}_n$ satisfies
%
%e5.2 #&#
\begin{equation}\label{eq:cond AN theta est}
0=\int\Psi(y,\widehat{\theta}_n)
P_n(\mathrm{d}y)+\mathrm{O}(n^{-1}),
\end{equation}
where $\Psi=(\Psi_1,\Psi_2,\Psi_3)$, defined as
%
%e5.3 #&#
\begin{eqnarray}
\label{eq:def Psi}
\Psi_1(y,\theta)&=& \mathbh{1}_{\{\|G^{-1}(y-m)\|\leq r\}
}G^{-1}(y-m),\nonumber\\
\Psi_2(y,\theta)&=& \mathbh{1}_{\{\|G^{-1}(y-m)\|\leq r\}}
\bigl(G^{-1}(y-m)(y-m)'G^{-1}-I_k \bigr),\\
\Psi_3(y,\theta)&=& \mathbh{1}_{\{\|G^{-1}(y-m)\|\leq r\}}-\gamma,\nonumber
\end{eqnarray}
where $\theta=(m,G,r)$, with $y,t\in\R^k$, $r>0$, and $G\in
\operatorname{PDS}(k)$. Rewrite equation \eqref{eq:cond AN theta est} as
%
%e5.4 #&#
\begin{eqnarray}\label{eq:expansion emp}
0&=&\Lambda(\widehat{\theta}_n) +
\int\Psi(y,\theta_0) (P_n-P_0)(\mathrm{d}y)\nonumber\\ [-8pt]\\ [-8pt]
&&{}+ \int\bigl(\Psi(y,\widehat{\theta}_n)-\Psi(y,\theta_0)\bigr)
(P_n-P_0)(\mathrm{d}y)+\mathrm{O}(n^{-1}),\nonumber
\end{eqnarray}
where
%
%e5.5 #&#
\begin{equation}\label{eq:def Lambda}
\Lambda(\theta)=\int\Psi(y,\theta) P_0(\mathrm{d}y).
\end{equation}
In order to determine the limiting distribution of
$\widehat{\theta}_n$, we proceed as follows. The first term on the
right-hand side of \eqref{eq:expansion emp} can be approximated by a
first order Taylor expansion that is linear in
$\widehat{\theta}_n-\theta_0$ and the second term can be treated by the
central limit theorem. Most of the difficulty is contained in the third
term, which must be shown to be of the order
$\mathrm{o}_{\mathbb{P}}(n^{-1/2})$. We apply empirical process theory, for
which we need $\int (\Psi(y,\widehat{\theta}_n)-\Psi(y,\theta_0) )
P_0(\mathrm{d}y)\to0$. For this, it suffices to impose
%
%e5.6 #&#
\begin{equation}\label{eq:mass E0}
P_0(\partial E_0)=0.
\end{equation}
For the MCD functional $\theta_n$ the argument is the same, apart from
the fact that replacing~$\phi_n$ by $\mathbh{1}_{E_n}$ requires an
additional condition on $P_0$, that is,
%
%e5.7 #&#
\begin{equation}\label{eq:atoms}
P_0\mbox{ has no atoms.}
\end{equation}
Note that \eqref{eq:mass E0} and \eqref{eq:atoms} are trivially true if
$P_0$ has a density. By representing elements of $\R^k\times
\operatorname{PDS}(k)\times\R$ as vectors, we then have the following central
limit theorem for the MCD estimators and the MCD functional at $P_n$.

%t5.1 #&#
\begin{theorem}
\label{th:AN} Let $P_0$ satisfy \eqref{eq:hyperprop}, \eqref{eq:cond
boundary} and \eqref{eq:mass E0}. Suppose that $ (\mu_0,\Sigma_0 )$ is
uniquely defined at~$P_0$. If $\Lambda$, as defined in \eqref{eq:def
Lambda}, has a non-singular derivative at $\theta_0$, then
\[
\widehat{\theta}_n-\theta_0
=-\Lambda'(\theta_0)^{-1}\frac1n\sum_{i=1}^n\bigl(\Psi(X_i,\theta
_0)-\mathbb{E}\Psi(X_i,\theta_0)\bigr)+\mathrm{o}_\mathbb{P}(n^{-1/2}),
\]
where $\Psi$ is defined in \eqref{eq:def Psi}. If in addition $P_0$
satisfies \eqref{eq:atoms}, then
\[
\theta_n-\theta_0
=-\Lambda'(\theta_0)^{-1}\frac1n\sum_{i=1}^n\bigl(\Psi(X_i,\theta
_0)-\mathbb{E}\Psi(X_i,\theta_0)\bigr)+\mathrm{o}_\mathbb{P}(n^{-1/2}).
\]
In particular, this means that $\sqrt{n}(\widehat{\theta}_n-\theta_0)$
and $\sqrt{n}(\theta_n-\theta_0)$ are asymptotically normal with mean
zero and covariance matrix
\[
\Lambda'(\theta_0)^{-1}M\Lambda'(\theta_0)^{-1},
\]
where $M$ is the covariance matrix of $\Psi(X_1,\theta_0)$.
\end{theorem}

Now Theorem \ref{th:AN} has been established, it turns out that the MCD
estimator and MCD functional (at $P_n$) are asymptotically equivalent,
that is, $\widehat{\theta}_n-\theta_n=\mathrm{o}_\mathbb{P}(n^{-1/2})$. Although
this seems natural, we have not been able to find an easier, direct
argument for this, in which case we could have avoided establishing
parallel results, such as the ones in Section \ref{subsec:consistency
MCD estimator}. An immediate consequence of Theorem \ref{th:AN} is
asymptotic normality of the MCD location estimator
$\sqrt{n}(\widehat{\mu}_n-\mu_0)$. Furthermore, since
\[
\widehat{\Sigma}_n-\Sigma_0=(\widehat{\Gamma}_n+\Gamma
_0)(\widehat{\Gamma}_n-\Gamma_0)=2\Gamma_0(\widehat{\Gamma
}_n-\Gamma_0)+\mathrm{o}_\mathbb{P}(1),
\]
Theorem \ref{th:AN} also yields asymptotic normality of the MCD
covariance estimator \mbox{$\sqrt{n}(\widehat{\Sigma}_n-\Sigma_0)$} and of
$\sqrt{n}(\widehat{\rho}_n-\rho_0)$. In \cite{catorlopuhaa2009}, a
precise expression is obtained for $\Lambda'(\theta_0)$ for $P_0$ with
a~density $f$ and nonsingularity of $\Lambda'(\theta_0)$ is proven if
$f$ has enough symmetry. This includes distributions with an
elliptically contoured density, so that as a special case of
Theorem~\ref{th:AN}, when $P_0$ has a unimodal elliptically contoured density,
one may recover Theorem 4 in \cite{butlerdaviesjuhn93} for the location
MCD estimator.

To determine the influence function, let $\phi_{\eps,x}$ be the
minimizing $\phi$-function for $P_{\eps,x}$ and let
\[\label{eq:def mu Gamma rho Pepsx}
\mu_{\eps,x}=T_{P_{\eps,x}}(\phi_{\eps,x}),\qquad
\Sigma_{\eps,x}=C_{P_{\eps,x}}(\phi_{\eps,x})=\Gamma_{\eps,x}^2,\qquad
\rho_{\eps,x}=r_{P_{\eps,x}}(\phi_{\eps,x})
\]
and $E_{\eps,x}=E(\mu_{\eps,x},\Sigma_{\eps,x},\rho_{\eps,x})$
be an
MCD functional at $P_{\eps,x}$ with corresponding minimizing ellipsoid.
To determine the influence function, we follow the same kind of
argument to obtain $M$-type score equations, by rewriting equations
\eqref{eq:def functionals} at $P_{\eps,x}$ and replacing
$\phi_{\eps,x}$ by $\mathbh{1}_{E_{\eps,x}}$. Note however, that from
the characterization given in Theorem \ref{th:characterization},
\[
\int(\mathbh{1}_{E_{\eps,x}}-\phi_{\eps,x})\,\mathrm{d}P_0 =
\cases{
P_0(\partial E_{\eps,x}), &\quad if $\phi_{\eps,x}=0$ on $\partial
E_{\eps,x}$,\cr
0, &\quad if $\phi_{\eps,x}=1$ on $\partial E_{\eps,x}$,\cr
P_0(\{z\}),&\quad otherwise, for some $z\in\partial E_{\eps,x}$.
}
\]
This means that in order to replace integrals over $\phi_{\eps,x}$ by
integrals over $E_{\eps,x}$, we need a~stronger condition on $P_0$,
that is,
%
%e5.8 #&#
\begin{equation}\label{eq:mass En}
P_0(\partial E)=0 \qquad\mbox{for any ellipsoid $E$}.
\end{equation}
Now, denote
$\theta_{\eps,x}=(\mu_{\eps,x},\Gamma_{\eps,x},\rho_{\eps,x})$ then
similar to \eqref{eq:cond AN theta est}, we obtain
%
%e5.9 #&#
\begin{equation}\label{eq:IF expansion theta}
0 = (1-\eps) \Lambda(\theta_{\eps ,x}) +
\eps\Phi_{\eps}(x),
\end{equation}
where $\Lambda$ is defined in \eqref{eq:def Lambda} and
$\Phi_{\eps}=(\Phi_{1,\eps},\Phi_{2,\eps},\Phi_{3,\eps})$, with
\begin{eqnarray*}
\Phi_{1,\eps}(x)&=& \phi_{\eps,x}(x)\Gamma_0^{-1}(x-\mu_0),\\
\Phi_{2,\eps}(x)&=& \phi_{\eps,x}(x) [\Gamma_0^{-1}(x-\mu
_0)(x-\mu_0)'\Gamma_0^{-1}-I_k ],\\
\Phi_{3,\eps}(x)&=& \phi_{\eps,x}(x)-\gamma.
\end{eqnarray*}
Define $\Theta(P)= (\mu(P),\Gamma(P),\rho(P) )$, where
$\mu(P)=T_P(\phi_P)$, $\Gamma(P)^2=C_P(\phi_P)$, $\rho(P)=r_P(\phi_p)$,
and $\phi_P$ denotes a minimizing trimming function. The influence
function of $\Theta(P)$ at $P_0$ is defined as
\[
\operatorname{IF}(x,\Theta,P_0) = \lim_{\eps\downarrow
0}\frac{\Theta((1-\eps)P_0+\eps\delta_x)-\Theta(P_0)}{\eps},
\]
if this limit exists, where $\delta_x$ is the Dirac measure at
$x\in\R^k$. The following theorem shows that this limit exists and
provides its expression.

%t5.2 #&#
\begin{theorem}
\label{th:IF} Suppose $P_0$ satisfies \eqref{eq:hyperprop},
\eqref{eq:cond boundary}, and \eqref{eq:mass En}. Suppose that $
(\mu_0,\Sigma_0 )$ is uniquely defined at $P_0$. Suppose that
$x\notin\partial E(\mu_0,\Sigma_0,\rho_0)$. If $\Lambda$ has a
non-singular derivative at $\theta_0$, then the influence function of
$\Theta$ at $P_0$ is given by
\[
\operatorname{IF}(x,\Theta,P_0)=-\Lambda'(\theta_0)^{-1}\Psi(x,\theta_0),
\]
where $\Psi$ is defined in \eqref{eq:def Psi}.
\end{theorem}

From definition \eqref{eq:def Psi}, we see that
$\operatorname{IF}(x,\Theta,P_0)$ is bounded uniformly for $x\notin
\partial
E(\mu_0,\Sigma_0,\rho_0)$. When $x\in\partial
E(\mu_0,\Sigma_0,\rho_0)$, then it is not clear what happens with
$\phi_{\eps,x}(x)$, as $\eps\downarrow0$. However, recall that there
exist $R>0$ such that $\{\phi_{\eps,x}\ne0\}\subset B_R$, for $\eps>0$
sufficiently small. This still implies that if $\phi_{\eps,x}(x)$ has a
limit, as $\eps\downarrow0$, then $\operatorname{IF}(x;\Theta,P_0)$ exists
and is bounded. In the case that $\phi_{\eps,x}(x)$ does not have a
limit, as $\eps\downarrow0$, then we can still conclude that
$\theta_{\eps,x}-\theta_0=\mathrm{O}(\eps)$, uniformly for $x\in\partial
E(\mu_0,\Sigma_0,\rho_0)$.

Because
$\Sigma_{\eps,x}-\Sigma_0=2\Gamma_0(\Gamma_{\eps,x}-\Gamma
_0)+o(1)$, as
$\eps\downarrow0$, it follows that the influence function of the
covariance functional $\Sigma(P)=C_P(\phi_P)$ is given by
\[
\operatorname{IF}(x;\Sigma,P_0) = 2\Gamma_0 \cdot\operatorname{IF}(x;\Gamma,P_0).
\]
As a special case of Theorem \ref{th:IF}, when $P_0$ has an
elliptically contoured density, Theorem~1 in \cite{crouxhaesbroeck99}
may be recovered (see \cite{catorlopuhaa2009}). Finally, note that
together with Theorem \ref{th:AN} it turns out that the von Mises
expansion indeed holds, that is,
\[\label{eq:von mises} \widehat{\theta}_n-\theta_0 = \frac1n \sum_{i=1}^n
\operatorname{IF}(X_i;\Theta,P_0)+ \mathrm{o}_\mathbb{P}(n^{-1/2}),
\]
which includes the heuristic \eqref{eq:von mises cov}.

On the basis of Theorems \ref{th:AN} and \ref{th:IF}, one could compute
asymptotic variances and robustness performance measures for the MCD
estimators and compare them with other robust competitors. Assuming the
influence function to exist and the expansion \eqref{eq:von mises cov}
to be valid, Croux and Haesbroeck \cite{crouxhaesbroeck99} provide an
extensive account of asymptotic and finite sample relative efficiencies
for the components of the MCD covariance estimator separately at the
multivariate standard normal, a contaminated multivariate normal and at
several multivariate Student distributions, for a variety of dimensions
$k=2,3,5,10,30$ and $\gamma=0.5,0.75$, as well as a comparison with
$S$-estimators and reweighted versions. Of particular interest would be
a comparison with the Stahel--Donoho (SD) estimator. Its asymptotic
properties have been established by Zuo \textit{et al.}
\cite{zuocuihe2004,zuocuiyoung2004,zuocui2005}, who also report an
asymptotic\vadjust{\goodbreak} and finite sample efficiency index for the SD location
estimator and for the full SD covariance estimator at the multivariate
normal and contaminated normal as well as a gross error sensitivity
index and maximum bias curve for the SD covariance estimator. The first
impression is that overall, apart from computational issues, the SD
estimator performs better than the MCD. However, a honest comparison
would require comparison of the same measure of efficiency and of the
maximum bias curves. To determine the latter seems far from trivial for
the MCD and we delay such a comparison to future research.

\begin{appendix}\label{sec:appendix}
%s6 #&#
\section*{Appendix}
Because the proof of Theorem
\ref{thm:contamination} relies heavily on the proof of the results in
Sections \ref{sec:existence} and \ref{sec:continuity}, this proof is
postponed to Section \ref{subsec:proofs continutuiy}.
\setcounter{equation}{0}

%s6.1 #&#
\subsection{\texorpdfstring{Proofs of existence and characterization (Section \protect\ref{sec:existence})}
{Proofs of existence and characterization (Section 3)}}
\label{subsec:proofs existence}

\begin{pf*}{Proof of Proposition \ref{prop:smalleig}}
Suppose that the
statement is not true. This means that the inequality either doesn't
hold for $t=0$ or for a sequence $t_n\to\infty$. Let us start with the
last case: we would have that there exist $a_n\in S^k$ and $\phi_n\in
K_{t_n}(\gamma)$ such that
%
%e6.1 #&#
\begin{equation}\label{eq:conv0}
\int\bigl(a_n'\bigl(x-T_{t_n}(\phi_n)\bigr)\bigr)^2\phi_n(x)P_n(\mathrm{d}x)
\to 0,
\end{equation}
where $P_n=P_{t_n}\to P_0$ weakly. Define the probability measures
\[
Q_n = \frac{\phi_n\cdot P_n}{\int\phi_n\,\mathrm{d}P_n}.
\]
Since $\{P_n\}$ is a tight (converging) sequence of probability
measures and since $\phi_n\leq1$ is such that $\int\phi_n\,\mathrm{d}P_n \geq
\gamma$ for all $n$, it is easy to see that $\{Q_n\}$ is a tight
sequence of probability measures. By going to subsequences, we may
assume that $Q_n\to Q_0$ weakly, for some probability measure $Q_0$.
Furthermore, again by going to subsequences, we may assume that $a_n\to
a_0\in S^k$ and that $\int\phi_n\,\mathrm{d}P_n\to\gamma_0\geq\gamma$.
Finally, using the tightness of $\{Q_n\}$, we can see that
$a_n'T_{t_n}(\phi_n)$ is a bounded sequence, since otherwise the limsup
of the integral would go to $+\infty$. This means that we can also
assume that $a_n'T_{t_n}(\phi_n)\to r\in\R$. Since $Q_n$ converges
weakly to $Q_0$, we can use the Skorokhod representation theorem to
define random variables $X_n$ and $X$ on a probability space $(\Omega,
\mu, \Sigma)$, such that $X_n\sim Q_n$, $X\sim Q_0$ and $X_n\to X$
almost surely. It follows that
\[
\bigl(a_n'\bigl(X_n - T_{t_n}(\phi_n)\bigr)\bigr)^2\to\bigl(a_0'(X-r)\bigr)^2,\qquad \mbox{a.s.},
\]
and from \eqref{eq:conv0}, together with Fatou's lemma we would
conclude that
\begin{eqnarray*}
\int\bigl(a_0'(x-r)\bigr)^2 Q_0(\mathrm{d}x) &=&
\mathbb{E}\bigl(a_0'(X-r)\bigr)^2\\
&\leq&\liminf_{n\to\infty}
\mathbb{E}\bigl(a_n'\bigl(X_n - T_{t_n}(\phi_n)\bigr)\bigr)^2\\
&=& \liminf_{n\to\infty} \frac1{\gamma_0}\int
\bigl(a_n'\bigl(x-T_{t_n}(\phi_n)\bigr)\bigr)^2\phi_n(x)P_n(\mathrm{d}x)=0.
\end{eqnarray*}
This would mean that $Q_0$ has to be concentrated on some closed
hyperplane $H$. However, we know that for all positive bounded
continuous functions $g$,
\[
\int g\,\mathrm{d}Q_0 = \lim_{n\to\infty} \int g\,\mathrm{d}Q_n \leq\lim_{n\to
\infty}\frac1{\gamma}\int g\,\mathrm{d}P_n = \frac1{\gamma} \int g\,\mathrm{d}P_0.
\]
This proves for the closed hyperplane $H$ that $1=Q_0(H)\leq
P_0(H)/\gamma$, which contradicts assumption \eqref{eq:hyperprop}. The
case $t=0$ can be proved using the same methods, just by defining $Q_n
=\phi_n\cdot P_0 / \int\phi_n\,\mathrm{d}P_0$.
\end{pf*}

The proof of Proposition \ref{prop:boundedsupp} relies on two
lemmas. The first one is a direct consequence of Proposition
\ref{prop:smalleig}, and shows that if $\det(C_t(\phi))$ is bounded
uniformly in $t$ and $\phi$, then there exists a fixed compact set that
contains all $(T_t(\phi),C_t(\phi))$ eventually. The second lemma is a
useful property involving the determinants of two nonnegative symmetric
matrices. Furthermore, for $R>0$ and $\mu\in\R^k$, define
%
%e6.2 #&#
\begin{equation}\label{eq:def B(mu,R)}
B(\mu,R) = \{x\in\R^k \dvt \|x-\mu\|\leq R\},
\end{equation}
and write $B_R$ in case $\mu=0$.

%l6.1 #&#
\begin{lemma}
\label{lem:compact set} Suppose $P_0$ satisfies \eqref{eq:hyperprop}
and let $P_t\to P_0$ weakly. Fix $M>0$. Then there exist $t_0\geq1$,
$0<\lambda_0\leq\lambda_1<\infty$ and $L,\rho>0$, such that for $t=0$,
all $t\geq t_0$, all $\phi$ such that~$T_t(\phi)$, $C_t(\phi)$ exist,
$\int\phi\,\mathrm{d}P_t\geq\gamma$, and
\[
\det(C_t(\phi))\leq M,
\]
we have that all eigenvalues of $C_t(\phi)$ are between $\lambda_0$ and
$\lambda_1$, $\|T_t(\phi)\|\leq L$, and $r_t(\phi)\leq\rho$.
\end{lemma}

\begin{pf}
The existence of $\lambda_0$ follows directly from
Proposition \ref{prop:smalleig}. This also implies that the largest
eigenvalue $\lambda_{\max}$ of $C_t(\phi)$ is smaller than
$M/\lambda_0^{k-1}$. Finally, choose $R>0$ such that for all $t\geq0$,
$P_t(B_R)\geq1-\gamma/2$, with $B_R$ as defined in \eqref{eq:def
B(mu,R)}. Suppose $\|T_t(\phi)\|\geq R$ and according to Lemma
\ref{lem:mass phi}, assume without loss of generality that $\int\phi\,\mathrm{d}P_t=\gamma$. Then, since
\[
\int_{B_R}\phi\,\mathrm{d}P_t = \int\phi\,\mathrm{d}P_t-\int_{\R^k\setminus B_R}\phi\,\mathrm{d}P_t
\geq\gamma-\bigl(1-P_t(B_R)\bigr) \geq\frac{\gamma}{2},
\]
we find
\begin{eqnarray*}
\label{eq:bound T}
\lambda_{\max} & \geq&
\frac1{\gamma}\int\biggl(\frac{T_t(\phi)'(T_t(\phi)-x)}{\|T_t(\phi)\|
}\biggr)^2\phi(x) P_t(\mathrm{d}x)\\
& \geq& \frac1{\gamma}\int_{B_{R}} \biggl(\frac{T_t(\phi)' (T_t(\phi)-x
)}{\|T_t(\phi)\|}\biggr)^2\phi(x) P_t(\mathrm{d}x)\\
& \geq&\frac1{\gamma}\int_{B_{R}} \bigl(\|T_t(\phi)\|-R\bigr)^2\phi(x) P_t(\mathrm{d}x)
\geq\frac{1}2 \bigl(\|T_t(\phi)\|-R\bigr)^2.
\end{eqnarray*}
This proves that there exists $L>0$, depending on $R$, $M$ and
$\lambda_0$, such that $\|T_t(\phi)\|\leq L$. Finally, since
\[
\bigl(x-T_t(\phi)\bigr)'C_t(\phi)^{-1}\bigl(x-T_t(\phi)\bigr)\leq
\frac{(\|x\|+L)^2}{\lambda_0}
\]
for $\rho>0$ large enough, the ellipsoid $E(T_t(\phi),C_t(\phi),\rho)$,
as defined in \eqref{eq:def ellipsoid}, contains the ball
$B(0,\rho\sqrt{\lambda_0}-L)$. Choose $\rho>0$ large enough, such that
\[
P_t\bigl(B\bigl(0,\rho\sqrt{\lambda_0}-L\bigr)\bigr)\geq\gamma
\]
for all $t\geq1$. Then $P_t(E(T_t(\phi),C_t(\phi),\rho))\geq\gamma$
and by definition we must have $r_t(\phi)\leq\rho$.
\end{pf}

%l6.2 #&#
\begin{lemma}\label{lem:mindet}
Let $\Sigma_1$ and $\Sigma_2$ be two symmetric matrices, nonnegative
and positive definite, respectively, such that
$\Tr(\Sigma_2^{-1}(\Sigma_1-\Sigma_2))<0$. Then
$\det(\Sigma_1)<\det(\Sigma_2)$. A similar result holds with $\leq$
instead of strict inequalities.
\end{lemma}

\begin{pf}
Without loss of generality we may assume that $\Sigma_2 =
I_k$. Suppose $\Tr(\Sigma_1)<\Tr(I_k)$. This means that the eigenvalues
$\lambda_1,\ldots,\lambda_k$ of $\Sigma_1$ satisfy
$(\lambda_1+\cdots+\lambda_k)/k < 1$. By means of the inequality
between the arithmetic mean and the geometric mean of nonnegative
numbers, we find $\det(\Sigma_1)=\lambda_1\cdots\lambda_k <1$.
\end{pf}

We also need the following well known result.
%
%l6.3 #&#
\begin{lemma}\label{lem:minmean}
Suppose $Q$ is a probability measure on $\R^k$ such that $\int\|x\|^2
Q(\mathrm{d}x) < +\infty$ and $Q$ is not supported by a hyperplane. Define $\mu
= \int x Q(\mathrm{d}x)$. Then for all $a\in\R^k$, $a\neq\mu$, we have
\[
\det\biggl(\int(x-\mu)(x-\mu)' Q(\mathrm{d}x)\biggr) < \det\biggl(\int(x-a)(x-a)' Q(\mathrm{d}x)\biggr).
\]
\end{lemma}

\begin{pf}
First note that
\[
\int(x-a)(x-a)' Q(\mathrm{d}x) = \int(x-\mu)(x-\mu)' Q(\mathrm{d}x) + (a-\mu)(a-\mu)',
\]
then apply Lemma \ref{lem:mindet}, remembering that $\int(x-a)(x-a)'
Q(\mathrm{d}x)$ is invertible and therefore strictly positive definite.
\end{pf}

\begin{pf*}{Proof of Proposition \ref{prop:boundedsupp}}
Choose
$t_0\geq1$ and $\lambda_0>0$ according to Proposition
\ref{prop:smalleig}. Choose $R'>0$ such that $P_t(B_{R'})\geq\gamma$,
for all $t\geq0$. Let $\psi_0$ be a continuous bounded function
such\vadjust{\goodbreak}
that $\mathbh{1}_{B_{R'}}\leq\psi_0\leq\mathbh{1}_{B_{R'+1}}$ and
define $D_0=2\det(C_0(\psi_0))$. Because $P_t\to P_0$ weakly, and
$\psi_0$ has bounded support, we have
\begin{eqnarray*}
&\displaystyle\int\psi_0\,\mathrm{d}P_t\to\int\psi_0\,\mathrm{d}P_0,\qquad
\int\psi_0(x)x P_t(\mathrm{d}x)\to\int\psi_0(x)x P_0(\mathrm{d}x)\quad \mbox{and}&\\[-3pt]
&\displaystyle\int\psi_0(x)xx' P_t(\mathrm{d}x)\to\int\psi_0(x)xx' P_0(\mathrm{d}x),&
\end{eqnarray*}
and hence $C_t(\psi_0)\to C_0(\psi_0)$, so that for $t$ large enough,
$\det(C_t(\psi_0))\leq D_0$.

Now, consider $\phi\in K_t(\gamma)$. If $\det(C_t(\phi))\geq
D_0\geq
\det(C_t(\psi_0))$, then we are done, because $\psi_0\in
K_t^{R'+1}(\gamma)$. Therefore, suppose that $\det(C_t(\phi))<D_0$.
According to Lemma \ref{lem:compact set}, this implies there exist
$\lambda_1\geq\lambda_0>0$ and $L>0$ such that
%
%e6.3 #&#
\begin{equation}\label{eq:bound lambda max and T}
\lambda_0 \leq\lambda_{\min}(C_t(\phi)) \leq\lambda_{\max}(C_t(\phi))\leq\lambda_1\quad
\mbox{and}\quad \|T_t(\phi)\|\leq L,
\end{equation}
uniformly in $t$ and $\phi$. According to Lemma \ref{lem:mass phi}, we
may assume that $\int\phi\,\mathrm{d}P_t=\gamma$. Choose any $R>R'+1$ and
suppose that $\phi>\mathbh{1}_{B_R}\cdot\phi$. Because $\int\phi\,\mathrm{d}P_t=\gamma$ and $P_t(B_{R'})\geq\gamma$, we know that
\[
\int_{B_{R'}}(1-\phi)\,\mathrm{d}P_t \geq\int_{\R^k\setminus B_R} \phi\,\mathrm{d}P_t.
\]
Define $h_1 = \mathbh{1}_{\R^k\setminus B_R}\cdot\phi$ and $h_2 =
\tau
\mathbh{1}_{B_{R'}}(1-\phi)$, where we choose $0\leq\tau\leq1$ such
that
%
%e6.4 #&#
\begin{equation}\label{eq:int h1=int h2}
\int h_2\,\mathrm{d}P_t = \int_{B_{R'}} \tau(1-\phi)\,\mathrm{d}P_t = \int h_1\,\mathrm{d}P_t.
\end{equation}
Furthermore, define $\psi= \phi- h_1 + h_2$ and note that $\psi\in
K^R_t(\gamma)$. Because according to~\eqref{eq:int h1=int h2}, $\int
\psi\,\mathrm{d}P_t=\int\phi\,\mathrm{d}P_t$, we can write
%
%e6.5 #&#
\begin{eqnarray}
\label{eq:psi phi argument}
\det(C_t(\psi)) & =&
\det\biggl(\frac1{\int\psi\,\mathrm{d}P_t} \int\bigl(x-T_t(\psi)\bigr)\bigl(x-T_t(\psi)\bigr)'\psi
(x) P_t(\mathrm{d}x)\biggr)\nonumber\\[-3pt]
& \leq&
\det\biggl(\frac1{\int\psi\,\mathrm{d}P_t} \int\bigl(x-T_t(\phi)\bigr)\bigl(x-T_t(\phi)\bigr)'\psi
(x) P_t(\mathrm{d}x)\biggr)\\[-3pt]
& =& \det\biggl(C_t(\phi) + \frac1{\int\phi\,\mathrm{d}P_t} \int
\bigl(x-T_t(\phi)\bigr)\bigl(x-T_t(\phi)\bigr)'\bigl(\psi(x)-\phi(x)\bigr) P_t(\mathrm{d}x)\biggr).\nonumber
\end{eqnarray}
For the inequality, we used Lemma \ref{lem:minmean}. So according to
Lemma \ref{lem:mindet}, it suffices to show that
%
%e6.6 #&#
\begin{equation}\label{eq:diff neg}
\frac1{\int\phi\,\mathrm{d}P_t} \int
\bigl(x-T_t(\phi)\bigr)'C_t(\phi)^{-1}\bigl(x-T_t(\phi)\bigr)\bigl(h_2(x)-h_1(x)\bigr) P_t(\mathrm{d}x) <0.
\end{equation}
To see that this is true, note that with \eqref{eq:bound lambda max and
T} we get
\begin{eqnarray*}
\int\bigl(x-T_t(\phi)\bigr)'C_t(\phi)^{-1}\bigl(x-T_t(\phi)\bigr)h_2(x) P_t(\mathrm{d}x) &
\leq&
{\lambda_0^{-1}}\int_{B_{R'}} \|x-T_t(\phi)\|^2h_2(x) P_t(\mathrm{d}x)\\[-3pt]
& \leq&{\lambda_0^{-1}}(R'+L)^2 \int h_2\,\mathrm{d}P_t,\vadjust{\goodbreak}
\end{eqnarray*}
and
\begin{eqnarray*}
\int\bigl(x-T_t(\phi)\bigr)'C_t(\phi)^{-1}\bigl(x-T_t(\phi)\bigr)h_1(x) P_t(\mathrm{d}x) &
\geq&
{\lambda_1^{-1}}\int_{\R^k\setminus B_{R}} \|x-T_t(\phi)\|^2h_1(x)
P_t(\mathrm{d}x)\\[-2pt]
& \geq&{\lambda_1^{-1}}(R-L)^2\int h_1\,\mathrm{d}P_t.
\end{eqnarray*}
So, together with \eqref{eq:int h1=int h2}, for $R$ large enough (but
independent of $\phi$!), this proves \eqref{eq:diff neg}.\vspace*{-2pt}
\end{pf*}

\begin{pf*}{Proof of Theorem \ref{th:existence}}
Choose $R>0$
and $t_0\geq1$ according to Proposition \ref{prop:boundedsupp}. Then
for $t=0$ and $t\geq t_0$, we may restrict minimization to $\phi\in
K_t^R(\gamma)$. Since $K^R_t(\gamma)$ is a weak*-compact subset of
$L^\infty(P_t|_{B_R})$, and since $\phi\mapsto\det(C_t(\phi))$ is a
weak*-continuous function on this space, we conclude that there exists
at least one minimum.\vspace*{-2pt}
\end{pf*}

\begin{pf*}{Proof of Theorem \ref{th:characterization}}
First, only
consider minimizing functions with $\int\phi\,\mathrm{d}P=\gamma$, which is
always possible according to Lemma \ref{lem:mass phi}. Let
$\phi_\gamma\in K_P(\gamma)$ be such that \mbox{$\int\phi_\gamma\,\mathrm{d}P=\gamma$}
and $\mathbh{1}_{E_P(\phi)^\circ}\leq\phi_\gamma\leq
\mathbh{1}_{E_P(\phi)}$. Define $Y=\Gamma_P(\phi)^{-1}(X-T_P(\phi))$,
where $X$ has distribution~$P$ and $\Gamma_P(\phi)^2=C_P(\phi)$, and
denote $Q$ for the distribution of $Y$. Define
%
%e6.7 #&#
\begin{eqnarray}\label{def:tau}
\tau_\gamma(y) &=&
\phi_\gamma\bigl(\Gamma_P(\phi)y+T_P(\phi)\bigr),\nonumber\\ [-9pt]\\ [-9pt]
\tau_\phi(y) &=& \phi\bigl(\Gamma_P(\phi)y+T_P(\phi)\bigr).\nonumber
\end{eqnarray}
Then $\mathbh{1}_{B^\circ}\leq\tau_\gamma\leq\mathbh{1}_{B}$, where
$B$ denotes the closed ball around the origin with radius~$r_P(\phi)$.
According to Lemma 1 in \cite{gordaliza1991}, with $\Phi(t)=t^2$, it
then follows that
\[
\int\tau_\gamma(y)\|y\|^2 Q(\mathrm{d}y) \leq\int\tau_\phi(y)\|y\|^2 Q(\mathrm{d}y)
\]
and one has equality if and only if $\mathbh{1}_{B^\circ}\leq
\tau_\phi\leq\mathbh{1}_{B}$ $Q$-a.s. This is equivalent with
\[
\int\bigl(\phi_\gamma(x)-\phi(x)\bigr) \bigl(x-T_P(\phi)\bigr)'C_P(\phi
)^{-1}\bigl(x-T_P(\phi)\bigr)
P(\mathrm{d}x) \leq0,
\]
where equality holds if and only if $\mathbh{1}_{E_P(\phi)^\circ
}\leq
\phi\leq\mathbh{1}_{E_P(\phi)}$, $P$-a.e. Then with Lemma
\ref{lem:mindet}, similar to arguments given in \eqref{eq:psi phi
argument} and \eqref{eq:diff neg}, it follows that
$\det(C_P(\phi_\gamma))<\det(C_P(\phi))$ unless
$\mathbh{1}_{E_P^\circ}\leq\phi\leq\mathbh{1}_{E_P}$, $P$-a.e.

Now, suppose that $\int\phi\,\mathrm{d}P>\gamma$. Then, according to Lemma
\ref{lem:mass phi}, for some $0<\lambda<1$, the function $0\leq
\lambda\phi<1$ would also be minimizing and satisfies
$\int(\lambda\phi)\,\mathrm{d}P=\gamma$. But then, the argument above shows that
$\lambda\phi=1$ on the interior of its own ellipsoid
$E(T_P(\lambda\phi),\allowbreak C_P(\lambda\phi),r_P(\lambda\phi))$, which is a
contradiction. We conclude that we must have $\int\phi\,\mathrm{d}P=\gamma$.

The last statement of the theorem is a little bit more subtle. Suppose
$P(\partial E_P)>0$, since otherwise the statement is trivially true.
Consider the following two functions on~$[0,1]$:
\[
f_1(t) = \frac{P(\partial E_P \cap\{\phi\leq t\})}{P(\partial E_P)}
\quad\mbox{and}\quad f_2(t)=\frac{P(\partial E_P \cap\{\phi\geq
t\})}{P(\partial E_P)}.
\]
If one realizes that $f_1+f_2\geq1$, $f_1$ is nondecreasing and
continuous from the right, whereas~$f_2$ is nonincreasing and\vadjust{\goodbreak}
continuous from the left, it is not hard to see that either $f_1=1$ on
$[0,1]$, in which case $\phi=0$, $P$-a.e. on $\partial E_P$, or $f_2=1$
on $[0,1]$, in which case $\phi=1$, $P$-a.e. on $\partial E_P$, or
there exists $t\in(0,1)$ such that $f_1(t),f_2(t)>0$. For this $t\in
(0,1)$, define
\[
A = \partial E_P \cap\{\phi\leq t\} \quad\mbox{and}\quad B = \partial E_P
\cap\{\phi\geq t\}.
\]
Either $P(A\cup B) = P(\{x\})$ for some $x\in\partial E_P$, in which
case $P(\partial E_P) = P(\{x\})$, or there exists $x\in\operatorname{supp}(P|_A)$
and $y\in\operatorname{supp}(P|_B)$ with $x\neq y$. We will show
that this last assumption will lead to a contradiction, thereby
finishing the proof. Choose $\eps>0$ such that $\eps< \|x-y\|/3$.
Define $A_\eps= A\cap B_\eps(x)$ and $B_\eps=B\cap B_\eps(y)$. By the
choice of $x$ and $y$, we know that $P(A_\eps),P(B_\eps)>0$. Choose
$\eta<\min(tP(B_\eps),(1-t)P(A_\eps))$ and define
\[
h_1(z) = \frac{\eta}{P(A_\eps)}\mathbh{1}_{A_\eps}(z)\quad
\mbox{and}\quad h_2(z) =
\frac{\eta}{P(B_\eps)}\mathbh{1}_{B_\eps}(z).
\]
Since $\phi\leq t$ on $A_\eps$ and $\phi\geq t$ on $B_\eps$, we get
that
\[
\psi= \phi+ h_1 - h_2 \in K_P(\gamma).
\]
Furthermore, $\int\psi\,\mathrm{d}P = \int\phi\,\mathrm{d}P = \gamma$. Since
$\eps<\|x-y\|/3$, we can see that
\[
T_P(\psi) = T_P(\phi) + \frac1{\gamma P(A_\eps)}\int_{A_\eps}z
P(\mathrm{d}z) -
\frac1{\gamma P(B_\eps)}\int_{B_\eps}z P(\mathrm{d}z) \neq T_P(\phi).
\]
Since $C_P(\psi)$ is invertible, this means that (with strict
inequality due to Lemma \ref{lem:minmean})
%
%e6.8 #&#
\begin{eqnarray}\label{eq:psi phi argument strict}
\det(C_P(\psi)) & =&
\det\biggl(\frac1{\gamma} \int\bigl(z-T_P(\psi)\bigr)\bigl(z-T_P(\psi)\bigr)'\psi(z)
P(\mathrm{d}z)\biggr)\nonumber\\
& <&
\det\biggl(\frac1{\gamma} \int\bigl(z-T_P(\phi)\bigr)\bigl(z-T_P(\phi)\bigr)'\psi(z)
P(\mathrm{d}z)\biggr)\\
& = & \det\biggl(C_P(\phi) + \frac1{\gamma} \int
\bigl(z-T_P(\phi)\bigr)\bigl(z-T_P(\phi)\bigr)'\bigl(h_1(z) - h_2(z)\bigr) P(\mathrm{d}z)\biggr).\nonumber
\end{eqnarray}
Since $A_\eps\cup B_\eps\subset\partial E_P$, and we know that for
$z\in\partial E_P$, $\operatorname{Tr}((z-T_P(\phi))'C_P(\phi)^{-1}(z-T_P(\phi)))$ is constant, we can use
Lemma \ref{lem:mindet} to conclude that
\[
\det(C_P(\psi)) < \det(C_P(\phi)),
\]
which contradicts the minimizing property of $\phi$.
\end{pf*}

%s6.2 #&#
\subsection{\texorpdfstring{Proofs of continuity (Section \protect\ref{sec:continuity})}
{Proofs of continuity (Section 4)}}
\label{subsec:proofs continutuiy} The proof of Theorem
\ref{th:continuity} uses the following two lemmas.
%
%l6.4 #&#
\begin{lemma}
\label{lem:cont1} Suppose $P_0$ satisfies \eqref{eq:hyperprop}. Let
$P_t\to P_0$ weakly and suppose that \eqref{eq:unif conv ellipsoid}
holds. For $t\geq1$, let $\psi_t\in K_t(\gamma)$ such that $\psi
_t\leq
\mathbh{1}_{E_t}$,\vadjust{\goodbreak} where $E_t=E(T_t(\psi_t),C_t(\psi_t),r_t(\psi_t))$,
and suppose there exists $R>0$, such that $\{\psi_t\ne0\}\subset B_R$,
for $t$ sufficiently large. Then\looseness=1
\[
\int\psi_t\,\mathrm{d}P_t- \int\psi_t\,\mathrm{d}P_0\to0,\qquad T_t(\psi_t)-T_0(\psi_t)\to0
\quad\mbox{and}\quad C_t(\psi_t)-C_0(\psi_t)\to0,
\]\looseness=0
uniformly in $\psi_t\in K_t(\gamma)$ with $\{\psi_t\ne0\}\subset B_R$.
\end{lemma}

\begin{pf}
Because $\{\psi_t\ne0\}\subset B_R$ eventually, we can
write
\[
\int\psi_t\,\mathrm{d}P_t- \int\psi_t\,\mathrm{d}P_0 = \int_{B_R}\psi_t\,\mathrm{d}(P_t-P_0)
\]
for $t$ sufficiently large. For the signed measure $Q_t=P_t-P_0$ write
$Q_t=Q_t^+-Q_t^-$, where~$Q_t^+$ and $Q_t^-$ are positive measures on
$\R^k$. According to \eqref{eq:unif conv ellipsoid},
\[
\sup_{E\in\mathcal{E}} Q_t^+(E)+\sup_{E\in\mathcal{E}} Q_t^-(E)
\leq
2 \sup_{E\in\mathcal{E}} |P_t(E)-P_0(E)| \to0.
\]
This implies $\sup_{E\in\mathcal{E}} Q_t^+(E)\to0$ and $\sup_{E\in
\mathcal{E}} Q_t^-(E)\to0$. Because $0\leq\psi_t\leq
\mathbh{1}_{E_t}$, we find
\begin{eqnarray*}
0&\leq&\int_{B_R}\psi_t(x) Q_t^+(\mathrm{d}x) \leq
Q_t^+(E_t\cap B_R)\leq Q_t^+(E_t)\to0,\\
0&\leq&\int_{B_R}\psi_t(x) Q_t^-(\mathrm{d}x) \leq Q_t^-(E_t\cap B_R)\leq
Q_t^-(E_t)\to0,
\end{eqnarray*}
which implies that
%
%e6.9 #&#
\begin{equation}\label{eq:int phi conv}
\int\psi_t\,\mathrm{d}P_t- \int\psi_t\,\mathrm{d}P_0 =
\int_{B_R}\psi_t(x) Q_t^+(\mathrm{d}x) - \int_{B_R}\psi_t(x) Q_t^-(\mathrm{d}x) \to0.
\end{equation}
Now, write
%
%e6.10 #&#
\begin{eqnarray}\label{eq:decompose diff T}
&&T_t(\psi_t)-T_0(\psi_t)\nonumber\\
&&\quad= \frac{1}{\int\psi_t\,\mathrm{d}P_t}\int_{B_R}\psi_t(x)x P_t(\mathrm{d}x) -
\frac{1}{\int\psi_t\,\mathrm{d}P_0}\int_{B_R}\psi_t(x)x P_0(\mathrm{d}x)\\
&&\quad= \frac{1}{\int\psi_t\,\mathrm{d}P_t}\int_{B_R}\psi_t(x)x (P_t-P_0)(\mathrm{d}x) + \biggl(
\frac{1}{\int\psi_t\,\mathrm{d}P_t}-\frac{1}{\int\psi_t\,\mathrm{d}P_0} \biggr)
\int_{B_R}\psi_t(x)x P_0(\mathrm{d}x).\nonumber
\end{eqnarray}
The first term in \eqref{eq:decompose diff T} tends to zero, because
$\gamma\leq\int\psi_t\,\mathrm{d}P_t\leq1$ and
\begin{eqnarray*}
0&\leq&\int_{B_R}\psi_t(x) \|x\| Q_t^+(\mathrm{d}x) \leq
RQ_t^+(E_t\cap B_R)\leq RQ_t^+(E_t)\to0,\\
0&\leq&\int_{B_R}\psi_t(x) \|x\| Q_t^-(\mathrm{d}x) \leq RQ_t^-(E_t\cap
B_R)\leq
RQ_t^-(E_t)\to0,
\end{eqnarray*}
which implies that
\[
\int_{B_R}\psi_t(x)x (P_t-P_0)(\mathrm{d}x) = \int_{B_R}\psi_t(x)x Q_t^+(\mathrm{d}x) -
\int_{B_R}\psi_t(x) x Q_t^-(\mathrm{d}x) \to0.
\]
The second term in \eqref{eq:decompose diff T} also tends to zero,
because of \eqref{eq:int phi conv} and the fact that
\[
\biggl\|\int_{B_R}\psi_t(x)x P_0(\mathrm{d}x)\biggr\|\leq R.
\]
It follows that $T_t(\psi_t)-T_0(\psi_t)\to0$. Similarly, one proves
$C_t(\psi_t)-C_0(\psi_t)\to0$.
\end{pf}

%
%l6.5 #&#
\begin{lemma}\label{lem:cont2}
Suppose $P_0$ satisfies \eqref{eq:hyperprop}. Let
$P_t\to P_0$ weakly and suppose that \eqref{eq:unif conv ellipsoid}
holds. For $t\geq1$, let $\psi_t\in K_t(\gamma)$ such that $\psi
_t\leq
\mathbh{1}_{E_t}$, where $E_t=E(T_t(\psi_t),C_t(\psi_t),r_t(\psi_t))$,
and suppose there exist $R>0$ such that $\{\psi_t\ne0\}\subset B_R$,
for $t$ sufficiently large. Then there exist a subsequence
$t_m\to\infty$ and $\psi^*\in K_0^R(\gamma)$, such that
\[
\lim_{m\to\infty}  (T_0(\psi_{t_m}),C_0(\psi_{t_m}) ) =
(T_0(\psi^*),C_0(\psi^*) ).
\]
\end{lemma}

\begin{pf}
Since $0\leq\psi_t\leq1$ and $\{\psi_t\ne0\}
\subset
B_R$, the $\psi_t$ can be viewed as elements of the class
\[
\mathfrak{L}_0^R=\{\psi\in L^{\infty}(P_0)\dvt 0\leq\psi\leq1,
\{\psi\ne0\}\subset B_R\},
\]
which is a weak*-compact subset of $L^\infty(P_0|_{B_R})$. Hence, there
exist a subsequence $(\psi_{t_m})$ that has a $\mathrm{weak}^*$ limit
in $\mathfrak{L}_0^R$, say $\psi^*$. This means that for any $g\in
L^1(P_0)$,
%
%e6.11 #&#
\begin{equation}\label{eq:weak limit}
\lim_{m\to\infty}\int\psi_{t_m}g\,\mathrm{d}P_0 =
\int\psi^* g\,\mathrm{d}P_0.
\end{equation}
In particular, $\int\psi_{t_m}\,\mathrm{d}P_0\to\int\psi^*\,\mathrm{d}P_0$. Because
$\int\psi_t\,\mathrm{d}P_t\geq\gamma$, together with Lemma \ref{lem:cont1} this
implies
\[
\int\psi^*\,\mathrm{d}P_0 = \lim_{m\to\infty} \int\psi_{t_m}\,\mathrm{d}P_0 \geq
\gamma
- \lim_{m\to\infty}\int\psi_{t_m}\,\mathrm{d}(P_{t_m}-P_0) = \gamma,
\]
so that $\psi^*\in K_0^R(\gamma)$. Finally, since the support of both
$\psi_{t_m}$ and $\psi^*$ lies in $B_R$, it follows from \eqref{eq:weak
limit} that
\begin{eqnarray*}
T_0(\psi_{t_m}) &=&
\frac1{\int\psi_{t_m}\,\mathrm{d}P_0}\int_{B_R} \psi_{t_m}(x)x P_0(\mathrm{d}x)\\
&=& \frac1{\int\psi^*\,\mathrm{d}P_0}\int_{B_R} \psi_{t_m}(x)x P_0(\mathrm{d}x) + \biggl(
\frac1{\int\psi_{t_m}\,\mathrm{d}P_0}-\frac1{\int\psi^*\,\mathrm{d}P_0} \biggr)
\int_{B_R} \psi_{t_m}(x)x P_0(\mathrm{d}x)\\
&\to&\frac1{\int\psi^*\,\mathrm{d}P_0} \int_{B_R} \psi^*(x)x P_0(\mathrm{d}x) =
T_0(\psi^*),
\end{eqnarray*}
and similarly $C_0(\psi_{t_m})\to C_0(\psi^*)$.\vadjust{\goodbreak}~%
\end{pf}

\begin{pf*}{Proof of Theorem \ref{th:continuity}}
Consider the
sequence $ (T_t(\phi_t),C_t(\phi_t) )$. According to Proposition
\ref{prop:smalleig} there exist $\lambda_0>0$, such that $\lambda _{\min}(C_t(\phi_t))\geq\lambda_0$ for $t$ sufficiently large. Similar to
the beginning of the proof of Proposition \ref{prop:boundedsupp}, we
obtain $\lambda_{\max}(C_t(\phi_t))\leq\lambda_1$ (see~\eqref{eq:bound lambda max and T}). Because $\psi_t\in K_t(\gamma)$,
again according to Proposition \ref{prop:smalleig}, $\lambda_{\min}(C_t(\psi_t))\geq\lambda_0$. Since $\det(C_t(\phi_t))\leq
\lambda_1^k$, for $t$ sufficiently large, and
$\det(C_t(\psi_t))-\det(C_t(\phi_t))$ tends to zero, it follows that
$\det(C_t(\psi_t))\leq2\lambda_1^k$ eventually, so that according to
Lemma \ref{lem:compact set}, there exists a compact set which contains
$ (T_t(\psi_t),C_t(\psi_t) )$ for $t$ sufficiently large. This means
there exist a convergent subsequence.

Now, consider a subsequence, which we continue to denote by $
(T_t(\psi_t),C_t(\psi_t) )$, for which $ (T_t(\psi_t),C_t(\psi_t)
)\to(T_0,C_0)$. From Lemmas \ref{lem:cont1} and \ref{lem:cont2}, we
conclude that there exists a further subsequence $(t_m)$, such that
\begin{eqnarray*}
T_0 &=& \lim_{m\to\infty}T_{t_m}(\psi_{t_m}) =
T_0(\psi^*),\\
C_0 &=& \lim_{m\to\infty} C_{t_m}(\psi_{t_m}) = C_0(\psi^*)
\end{eqnarray*}
for some $\psi^*\in K_0^R(\gamma)$. It remains to show that $(T_0,C_0)$
is an MCD-functional, that is, $\det(C_0(\psi^*))$ minimizes
$\det(C_0(\phi))$ over $K_0^R(\gamma)$. To this end, suppose
there exists $\delta>0$ and $\phi\in K_0^R(\gamma)$, such that
\[
\det(C_0(\phi))\leq\det(C_0(\psi^*))-\delta.
\]
Since the set of bounded continuous functions is dense within
$K_0^R(\gamma)$, we can construct a bounded continuous function
${\psi}\in K_0^R(\gamma)$, such that for all $i,j=1,2,\ldots,k$:
\[
\int|\psi-\phi|\,\mathrm{d}P_0,\qquad \int|\psi(x)x_i-\phi(x)x_i| P_0(x)
\quad\mbox{and}\quad \int|\psi(x)x_ix_j-\phi(x)x_ix_j| P_0(x),
\]
can be made arbitrarily small. Hence, we can construct a bounded
continuous function ${\psi}\in K_0^R(\gamma)$, such that
\[
\det(C_0({\psi}))\leq\det(C_0(\psi^*))-\delta/2.
\]
Now, since ${\psi}(x)x$ is bounded and continuous on $B_R$, we have
\[
T_{t_m}({\psi}) = \frac1{\int\psi\,\mathrm{d}P_{t_m}}\int{\psi}(x)x P_{t_m}(\mathrm{d}x)
\to\frac1{\int\psi\,\mathrm{d}P_{0}}\int{\psi}(x)x P_0(\mathrm{d}x) = T_0({\psi}),
\]
and similarly $C_{t_m}({\psi})\to C_0({\psi})$. Since also
$\det(C_{t_m}(\psi_{t_m}))-\det(C_{t_m}(\phi_{t_m}))\to0$, it would
follow that
\begin{eqnarray*}
\lim_{m\to\infty} \det(C_{t_m}({\psi})) &=&
\det(C_0({\psi})) \leq
\det(C_0(\psi^*))-\frac{\delta}2\\
&=& \lim_{m\to\infty}\det(C_{t_m}(\psi_{t_m}))-\frac
{\delta}2 =
\lim_{m\to\infty}\det(C_{t_m}(\phi_{t_m}))-\frac{\delta}2.
\end{eqnarray*}
This would mean that for $m$ sufficiently large,
$\det(C_{t_m}({\psi}))\leq
\det(C_{t_m}(\phi_{t_m}))-\delta/4$, which contradicts the
minimizing property of $\phi_{t_m}$.
\end{pf*}

\begin{pf*}{Proof of Theorem \ref{thm:contamination}}
First note
that when $\eps\downarrow0$, then $P_{r,\eps}\to P$ weakly. Condition~\eqref{eq:unif conv ellipsoid} automatically holds, and because
$P(H)<(\gamma-\eps)/(1-\eps)<\gamma$, also condition~\eqref{eq:hyperprop} holds. According to Theorem \ref{th:existence}
this means that $\operatorname{MCD}_\gamma(P_{r,\eps})$ exists, for $\eps>0$
sufficiently small, and the minimizing $\phi_{r,\eps}\in
K_{P_{r,\eps}}^R(\gamma)$. Hence, together with Theorem
\ref{th:characterization}, all conditions of Theorem
\ref{th:continuity} are satisfied, which yields the first limit in (i).
The proof the second limit in (i) mimics the proof of Theorem
\ref{th:continuity}. Note that although we are not dealing with a
weakly convergent sequence of measures satisfying condition
\eqref{eq:unif conv ellipsoid}, we do have that for all continuous
functions $f$ with bounded support,
%
%e6.12 #&#
\begin{equation}\label{eq:conv bounded f}
\lim_{\|r\|\to\infty} \int f\,\mathrm{d}P_{r,\eps} =
\int f\,\mathrm{d}(1-\eps)P
\end{equation}
and for every fixed $R>0$,
%
%e6.13 #&#
\begin{equation}\label{eq:unif conv ellipoid BR}
\lim_{\|r\|\to\infty} \sup_{E\in
\mathcal{E}, E\subset B_R} |P_{r,\eps}(E) - (1-\eps)P(E)| = 0.
\end{equation}
We first show the analogue of Proposition \ref{prop:boundedsupp}, that
is, there exists $R>0$ such that for all $\|r\|$ sufficiently large,
the support of all minimizing $\phi$ for $P_{r,\eps}$ lies in $B_R$.

Choose $R'>0$ large enough, such that $P(B_{R'})>\gamma/(1-\eps)$. This
shows that there exist $\psi\in K_{\gamma/(1-\eps)}(P)$ with support
contained in $B_{R'}$, such that $\int\psi\,\mathrm{d}P_{r,\eps} \geq\gamma$,
for all $r\in\R^k$, and from \eqref{eq:conv bounded f} we find
\[
\lim_{\|r\|\to\infty} (T_{P_{r,\eps}}(\psi),C_{P_{r,\eps}}(\psi
)) =
\bigl(T_{(1-\eps)P}(\psi),C_{(1-\eps)P}(\psi)\bigr).
\]
When we take $M=2\det(C_{(1-\eps)P}(\psi))$, there exists
$r_0>0$, such that for all $r$ with \mbox{$\|r\|>r_0$},
$\det(C_{P_{r,\eps}}(\psi))<M$. It follows, that if $\phi$ is a
minimizing function for $P_{r,\eps}$ at level $\gamma$, we can conclude
that $\det(C_{P_{r,\eps}}(\phi))<M$, for $\|r\|>r_0$. Also, since
$\int
\phi\,\mathrm{d}P \geq(\gamma- \eps)/(1-\eps)$, Proposition~\ref{prop:smalleig}
yields that there exists $\lambda_0>0$, not depending on $\phi$, such
that for all $a\in\R^k$ with $\|a\|=1$
\[
\int\bigl(a'\bigl(x-T_P(\phi)\bigr)\bigr)^2\phi\,\mathrm{d}P \geq\lambda_0.
\]
This implies that
\begin{eqnarray*}
\int\bigl(a'\bigl(x-T_{P_{r,\eps}}(\phi)\bigr)\bigr)^2\phi\,\mathrm{d}P_{r,\eps} &
\geq&
(1-\eps) \int\bigl(a'\bigl(x-T_{P_{r,\eps}}(\phi)\bigr)\bigr)^2\phi\,\mathrm{d}P\\
& \geq& (1-\eps) \int\bigl(a'\bigl(x-T_{P}(\phi)\bigr)\bigr)^2\phi\,\mathrm{d}P \geq
(1-\eps)\lambda_0.
\end{eqnarray*}
This means that for all minimizing $\phi$ for $P_{r,\eps}$, we have a
uniform lower bound on the smallest eigenvalue. From here on, we copy
the proof of Lemma \ref{lem:compact set}. We choose \mbox{$R>0$} and
$\delta>0$ (independent of $\phi$) such that $P_{r,\eps}(B_R)\geq
(1-\eps)P(B_R)\geq1-\gamma+\delta$ and \mbox{$\int_{B_R}\phi\,\mathrm{d}P_{r,\eps
}\geq\delta$}. This then shows that there exists $\lambda_{\max}$ and
$L>0$ such that for all $r$ with $\|r\|>r_0$ and for all minimizing
$\phi$, we have $\|T_{P_{r,\eps}}(\phi)\|<L$ and the largest eigenvalue
of $C_{P_{r,\eps}}(\phi)$ is smaller than\vadjust{\goodbreak} $\lambda_{\max}$. Now we
can follow the proof of Proposition~\ref{prop:boundedsupp}, starting
from \eqref{eq:bound lambda max and T}, to conclude that there exists
$R>0$ and $r_0>0$, such that for all $r$ with $\|r\|>r_0$, the support
of all minimizing $\phi$ for $P_{r,\eps}$ lies within $B_R$.

Note that, because according to Proposition \ref{prop:boundedsupp} all
minimizing $\phi$ for $(1-\eps)P$ also have a~fixed bounded support,
this immediately yields statement (ii). Indeed, if $Q$ has bounded
support, then for $\|r\|$ sufficiently large, $\int\phi\,\mathrm{d}Q\circ
\tau_r^{-1}=0$ for all $\phi$ with a fixed bounded support. Hence, for
$\|r\|$ sufficiently large, $\phi$ is minimizing for $P_{r,\eps}$ at
level $\gamma$ if and only if $\phi$ is minimizing for $(1-\eps)P$ at
level $\gamma$, which means that $\phi$ is minimizing for~$P$ at level
$\gamma/(1-\eps)$.

To finish the proof of (i), we follow the proof of Theorem
\ref{th:continuity}, from the point of considering a~convergent
subsequence of $\operatorname{MCD}_\gamma(P_{r,\eps})$. The conclusions of
Lemmas \ref{lem:cont1} and \ref{lem:cont2} are still valid if we
replace the condition of weak convergence by \eqref{eq:conv bounded f}
and replace condition \eqref{eq:unif conv ellipsoid} by \eqref{eq:unif
conv ellipoid BR}. This means that the proof of the second limit in (i)
is completely similar to the remainder of the proof of Theorem
\ref{th:continuity}, which proves (i).
\end{pf*}

\begin{pf*}{Proof of Corollary \ref{cor:consistency functional}}
Since the MCD functional at $P_0$ is unique, it follows immediately
from Theorem \ref{th:continuity} that each convergent subsequence has
the same limit point $(T_0(\phi_0),C_0(\phi_0))$, which proves part
(i).

For $t=1,2,\ldots,$ write $E_{t}(s)=E(T_t(\psi_t),C_t(\psi_t),s)$ and
$\rho_t=r_t(\psi_t)$, as defined by \eqref{eq:def ellipsoid} and
\eqref{eq:def r_n}, and write $E_0(s)$ and $\rho_0$ for the ellipsoid
and radius corresponding to $\phi_0$. Let $\eps>0$ and suppose that
$\liminf\rho_t\leq\rho_0-3\eps$. Then there exists a subsequence
$t_n\to\infty$ such that\looseness=-1
\[
E_{t_n}(\rho_{t_n})\subset E_{t_n}(\rho_0-2\eps) \subset
E_0(\rho_0-\eps),
\]\looseness=0
and from \eqref{eq:cond boundary} it would follow that
$P_{t_n}(E_{t_n}(\rho_{t_n}))\leq P_{t_n}(E_0(\rho_0-\eps))\to
P_0(E_0(\rho_0-\eps))<\gamma$, which is in contradiction with the
definition of $\rho_{t_n}$. Similarly, if $\limsup\rho_{t}\geq
\rho_0+3\eps$, then there exists a subsequence $t_n\to\infty$ such that
\[
E_{t_n}(\rho_{t_n}-\eps)\supset E_{t_n}(\rho_0+\eps) \supset
E_0(\rho_0),
\]
from which it would follow that $P_{t_n}(E_{t_n}(\rho_{t_n}-\eps
))\geq
P_{t_n}(E_{0}(\rho_0))\to P_{0}(E_{0}(\rho_0))\geq\gamma$, which again
is in contradiction with the definition of $\rho_{t_n}$. Since $\eps>0$
was arbitrary, this finishes the proof of part (ii).
\end{pf*}

\begin{pf*}{Proof of Proposition \ref{prop:diff det
subsample-phi}}
Let $\phi_n$ be a minimizing function for the MCD
functional corresponding to $P_n$. Then by definition
\[
\det(C_n(\phi_n)) \leq\det(C_n(\mathbh{1}_{S_n})) =
\det(\widehat{C}_n(S_n)).
\]
First note that $\phi_n$ cannot be zero on the boundary of
$E_n=E(T_n(\phi_n),C_n(\phi_n),r_n(\phi_n))$. Hence according to
Theorem \ref{th:characterization}, we either have $\phi_n=1$ on
$\partial E_n$ or there exists a point $x\in\partial E_n$ such that
$P_n(\{x\})=P_n(\partial E_n)$. In the first case
\[
\det(\widehat{C}_n(S_n)) \leq\det(\widehat{C}_n(\widetilde{S}_n)) =
\det(C_n(\phi_n))
\]
for the subsample $\widetilde{S}_n=\{X_i\dvt\phi_n(X_i)=1\}$, which means
$\det(\widehat{C}_n(S_n))=\det(C_n(\phi_n))$.

Consider the other case. Suppose $\phi_n=1$ in $k$ points other than
$x$ and suppose there are~$m$ sample\vadjust{\goodbreak} points $X_i=x$. Then we must have
$\gamma>k/n$ and $\phi_n(x)=\eps_n$ for some $0<\eps_n<1$, where
\[
\gamma=\int\phi_n\,\mathrm{d}P_n=\frac{k}{n}+\frac{m\eps_n}n.
\]
Now, let $\widetilde{S}_n$ be the subsample consisting of the $k$
points where $\phi_n=1$ and $\lceil m\eps_n\rceil$ points $X_i=x$. Then
$\widetilde{S}_n$\vspace*{1pt} has $nP_n(\widetilde{S}_n)=k+\lceil
m\eps_n\rceil=\lceil n\gamma\rceil$ points. According to Proposition
\ref{prop:boundedsupp}, with probability one, there exists $R>0$ such
that $\widetilde{S}_n$ and $\{\phi_n\ne0\}$ are contained in $B_R$.
This implies
\[
\|T_n(\mathbh{1}_{\widetilde{S}_n})-T_n(\phi_n)\|=\frac{nR}{\lceil
n\gamma\rceil}\biggl(\frac{\lceil m\eps_n\rceil}{n}-\frac{m\eps_n}{n}\biggr)
\leq
\frac{R}{\lceil n\gamma\rceil},
\]
and similarly
$C_n(\mathbh{1}_{\widetilde{S}_n})-C_n(\phi_n)=\mathrm{O}(n^{-1})$ and
$\det(C_n(\mathbh{1}_{\widetilde{S}_n}))-\det(C_n(\phi_n))=\mathrm{O}(n^{-1})$,
with probability one. This means
\[
\det(\widehat{C}_n(S_n)) \leq\det(\widehat{C}_n(\widetilde{S}_n)) =
\det(C_n(\mathbh{1}_{\widetilde{S}_n})) = \det(C_n(\phi_n))+\mathrm{O}(n^{-1})
\]
with probability one.
\end{pf*}

The proof of Proposition \ref{prop:characterization
estimator} relies partly on the following property.
%
%l6.6 #&#
\begin{lemma}
\label{lem:size subsample} Let $S_m$ be a subsample of size $m\geq2$
and let $X^*\in S_m$ have maximal Mahalanobis distance with respect to
the corresponding trimmed sample mean $T_m=\widehat{T}_n(S_m)$ and
trimmed sample covariance $C_m=\widehat{C}_n(S_m)$, that is,
\[
X^*=\argmax\limits_{X_i\in S_m} (X_i-T_m)'C_m^{-1}(X_i-T_m).
\]
Define subsample $S_{m-1}=S_m\setminus\{X^*\}$ with trimmed sample
covariance $C_{m-1}=\widehat{C}_n(S_{m-1})$. Then
\[
\det(C_{m-1})\leq\det(C_m).
\]
\end{lemma}

\begin{pf}
We can write
\begin{eqnarray*}
\det(C_{m-1}) &=& \det\biggl( \frac1{m-1}\sum_{X_i\in
S_{m-1}}(X_i-T_{m-1})(X_i-T_{m-1})'\biggr)\\
&<& \det\biggl( \frac1{m-1}\sum_{X_i\in S_{m-1}}(X_i-T_{m})(X_i-T_{m})'
\biggr)\\
&=& \det\biggl( \frac{m}{m-1}C_m-\frac1{m-1}(X^*-T_m)(X^*-T_m)'\biggr)\\
&=& \det\biggl( C_m+\frac1{m-1} [C_m-(X^*-T_m)(X^*-T_m)' ] \biggr).
\end{eqnarray*}
From the definition of $C_m$, after multiplication with $C_m^{-1}$ and
taking traces, we find
\[
k=\frac1m\sum_{X_i\in S_m}(X_i-T_m)'C_m^{-1}(X_i-T_m).
\]
Therefore, since $X^*$ has the largest value for
$(X_i-T_m)'C_m^{-1}(X_i-T_m)$, it follows that
\[
\operatorname{Tr}[I_k-C_m^{-1}(X^*-T_m)(X^*-T_m)'] =
k-(X^*-T_m)'C_m^{-1}(X^*-T_m) \leq 0.
\]
The lemma now follows from Lemma \ref{lem:mindet}.
\end{pf}

\begin{pf*}{Proof of Proposition \ref{prop:characterization
estimator}}
Suppose that there is a point $X_\ell\in\widehat{E}_n$
that is not in $S_n$. Because $S_n$ must always have at least one point
on the boundary of $\widehat{E}_n$, we can then interchange a point
$X_j\in S_n$ that lies on the boundary of $\widehat{E}_n$ with
$X_\ell$. We will show that this will always decrease
$\det(\widehat{C}_n(S_n))$. Let
$S_n^*=(S_n\setminus\{X_j\})\cup\{X_\ell\}$. Then
\[
\widehat{T}_n(S_n^*) = \frac{1}{\lceil n\gamma\rceil} \sum_{X_i\in
S_n^*} X_i = \widehat{T}_n(S_n)+\frac{1}{\lceil
n\gamma\rceil}(X_j-X_\ell) \ne\widehat{T}_n(S_n).
\]
Therefore (with a strict inequality), similar to \eqref{eq:psi phi
argument strict}, we have
\begin{eqnarray*}
\det(\widehat{C}_n(S_n^*)) &<& \det\biggl( \widehat{C}_n(S_n) -
\frac{1}{\lceil n\gamma\rceil}
\bigl(X_j-\widehat{T}_n(S_n)\bigr)\bigl(X_j-\widehat{T}_n(S_n)\bigr)'\\
&&\hphantom{\det\biggl(\biggr.}{}+ \frac{1}{\lceil n\gamma\rceil}
\bigl(X_\ell-\widehat{T}_n(S_n)\bigr)\bigl(X_\ell-\widehat{T}_n(S_n)\bigr)' \biggr).
\end{eqnarray*}
Because $X_j$ is on the boundary of $\widehat{E}_n$ and $X_\ell$ inside
$\widehat{E}_n$, we have
\[
\bigl(X_\ell-\widehat{T}_n(S_n)\bigr)'\widehat{C}_n(S_n)^{-1}\bigl(X_\ell-\widehat
{T}_n(S_n)\bigr)
-
\bigl(X_j-\widehat{T}_n(S_n)\bigr)'\widehat{C}_n(S_n)^{-1}\bigl(X_j-\widehat{T}_n(S_n)\bigr)
\leq0.
\]
Therefore, it follows from Lemma \ref{lem:mindet} that
$\det(\widehat{C}_n(S_n^*))<\det(\widehat{C}_n(S_n))$, which
contradicts the minimizing property of $S_n$. We conclude that
$\{X_1,\ldots, X_n\}\cap\widehat{E}_n\subset S_n$. Since according to
Lemma \ref{lem:size subsample} the subsample $S_n$ has exactly $\lceil
n\gamma\rceil$ points, and by definition~$\widehat{E}_n$ contains at
least $\lceil n\gamma\rceil$ points, we conclude that $\{X_1,\ldots,
X_n\}\cap\widehat{E}_n = S_n$.
\end{pf*}

%l6.7 #&#
\begin{lemma}
\label{lem:bounded supp estimator} Suppose $P_0$ satisfies
\eqref{eq:hyperprop}. With probability one, there exists $R>0$ and
$n_0\geq1$, such that for all $n\geq n_0$ and all
subsamples $S_n$ with at least $n\gamma$ points, there exists a
subsample $S_n^*$ with exactly $\lceil n\gamma\rceil$ points contained
in $B_R$ such that
\[
\det(\widehat{C}_n(S_n^*))\leq\det(\widehat{C}_n(S_n)).
\]
\end{lemma}

\begin{pf}
The proof is along the lines of the proof of
Proposition \ref{prop:boundedsupp}. We first choose $R'>0$ and
construct a subsample $S_{n0}\subset B_{R'}$ with at least $n\gamma$
points, for which $\det(\widehat{C}_n(S_{n0}))$ is uniformly
bounded\vadjust{\goodbreak}
for $n$ sufficiently large. By the law of large numbers, $P_n\to P_0$
weakly with probability one. Hence, we can choose $R'>0$ such that for
all $n\geq1$,
\[
P_n(B_{R'})\geq\max\{1-\gamma/2,\gamma+(1-\gamma)/2\},
\]
with probability one, and define subsample $S_{n0}=\{X_i\dvt X_i\in
B_{R'}\}$. Then, $S_{n0}\subset B_{R'}$ and because $P_n(B_{R'})\geq
\gamma$, we also have $\mathbh{1}_{S_{n0}}\in K_n(\gamma)$. Then,
according to Proposition \ref{prop:smalleig}, with probability one,
there exist a $\lambda_0>0$ such that
\[
\lambda_{\min}(C_n(\mathbh{1}_{S_{n0}}))\geq\lambda_0
\]
for $n$ sufficiently large. Define
$D_0=2\det(C_0(\mathbh{1}_{B_{R'}}))$. From \eqref{eq:unif conv
ellipsoid}, we have $P_n(B_{R'})\to P_0(B_{R'})$, with probability one,
and since the functions $x$ and $xx'$ bounded and continuous on
$B_{R'}$, we also have
\[
\int_{B_{R'}}x\,\mathrm{d}P_n\to\int_{B_{R'}}x\,\mathrm{d}P_0,\quad \mbox{and}\quad
\int_{B_{R'}}xx'\,\mathrm{d}P_n\to\int_{B_{R'}}xx'\,\mathrm{d}P_0,
\]
with probability one. Hence, together with \eqref{eq:rel est-func}, it
follows that for $n$ sufficiently large,
\[
\det(\widehat{C}_n(S_{n0}))=\det(C_n(\mathbh{1}_{S_{n0}}))=\det
(C_n(\mathbh{1}_{B_{R'}}))\leq
D_0,
\]
with probability one. Now, let $S_n$ be a subsample with $h_n\geq
n\gamma$ points. According to Lemma \ref{lem:size subsample}, without
loss of generality, we may assume that is has exactly $\lceil
n\gamma\rceil$ points. When $\det(\widehat{C}_n(S_n))>D_0$, then we are
done because the subsample $S_{n0}$ has a smaller determinant, it is
contained in $B_{R'}$, and according to Lemma \ref{lem:size subsample}
we can reduce $S_{n0}$ if necessary to have exactly $\lceil
n\gamma\rceil$ points, without increasing the determinant. So suppose
that $S_n$ has $\lceil n\gamma\rceil$ points and
$\det(\widehat{C}_n(S_n))\leq D_0$. From here on the proof is identical
to that of Proposition \ref{prop:boundedsupp} and is left to the
reader.
\end{pf}

\begin{pf*}{Proof of Theorem \ref{th:consistency estimator}}
With probability one $P_n\to P_0$ weakly and \eqref{eq:unif conv
ellipsoid} holds, since the class of ellipsoids has polynomial
discrimination or forms a Vapnik--Cervonenkis class. According to
\eqref{eq:rel est-func} the MCD estimators can be written as MCD
functionals with trimming function $\psi_n=\mathbh{1}_{S_n}$. From
Propositions \ref{prop:diff det subsample-phi} and
\ref{prop:characterization estimator} together with Lemma
\ref{lem:bounded supp estimator}, it follows that $\psi_n$ satisfies
the conditions of Theorem \ref{th:continuity} with probability one,
which proves the theorem.
\end{pf*}

%s6.3 #&#
\subsection{\texorpdfstring{Proofs of asymptotic normality and IF (Section \protect\ref{sec:asymp norm})}
{Proofs of asymptotic normality and IF (Section 5)}}
\label{subsec:proofs AN} The proof of Theorem \ref{th:AN} relies on the
following result from \cite{pollard84}, which we state for easy
reference.
%
%t6.1 #&#
\begin{theorem}[(Pollard, 1984)]\label{th:pollard}Let $\mathcal{F}$ be a
permissible class of real valued functions with envelope $H\geq
|\phi|$, $\phi\in\mathcal{F}$, and suppose that
$0<\mathbb{E}[H(X)^2]<\infty$. If the class of graphs of functions
in\vadjust{\goodbreak}
$\mathcal{F}$ has polynomial discrimination, then for each $\eta>0$ and
$\eps>0$ there exists $\delta>0$ such that
\[
\limsup_{n\to\infty} \mathbb{P} \biggl\{ \sup_{(\phi_1,\phi_2)\in
[\delta]}
n^{1/2} \biggl| \int(\phi_1-\phi_2)\,\mathrm{d}(P_n-P_0) \biggr|
>\eta
\biggr\} <\eps,
\]
where $[\delta]=\{(\phi_1,\phi_2)\dvt\phi_1,\phi_2\in\mathcal
{F}\mbox{ and }\int(\phi_1-\phi_2)^2\,\mathrm{d}P_0<\delta^2\}$.
\end{theorem}
The theorem is not stated as such in \cite{pollard84}, but it is a
combination of the Approximation lemma (page 27), Lemma II.36 (page 34) and
the Equicontinuity lemma (page 150). The polynomial discrimination of
$\mathcal{F}$ provides a suitable bound on the entropy of $\mathcal{F}$
(Approximation lemma together with Lemma II.36). The stochastic
equicontinuity stated in Theorem \ref{th:pollard} is then a consequence
of the fact that the entropy of $\mathcal{F}$ is small enough
(Equicontinuity lemma). The classes of functions we will encounter in
this way can always be indexed by the parameter set $\R^k\times
\operatorname{PDS}(k)\times\R_+$, and are easily seen to be permissible in
the sense of Pollard \cite{pollard84}.

\begin{pf*}{Proof of Theorem \ref{th:AN}}
Consider equation
\eqref{eq:expansion emp} and define
\[
\mathcal{F}= \bigl\{ \mathbh{1}_{\{\|G^{-1}(x-m)\|\leq r\}}\dvt m\in\R^k,
G\in
\operatorname{PDS}(k), r>0 \bigr\}.
\]
As subclass of the class of indicator functions of all ellipsoids, the
class of graphs $\mathcal{G}$ of functions in $\mathcal{F}$ has
polynomial discrimination and obviously $\mathcal{F}$ has envelope
$H=1$. Hence, Theorem \ref{th:pollard} applies to $\Psi_3$. For the
real valued components of $\Psi_1$ and $\Psi_2$, use that there exists
$R>0$, such that for $n=0$ and $n$ sufficiently large
\[
\{x\in\R^k\dvt\|\Gamma_n^{-1}(x-\mu_n)\|\leq\rho_n\}\subset B_R.
\]
This means that for all $i,j=1,2,\ldots,k$, the classes
\begin{eqnarray*}
\mathcal{F}_i &=& \bigl\{ x_i\mathbh{1}_{\{\|G^{-1}(x-m)\|\leq r\}\cap
B_R}\dvt
m\in\R^k, G\in\operatorname{PDS}(k), r>0\bigr\},\\[-2pt]
\mathcal{F}_{ij} &=& \bigl\{ x_ix_j\mathbh{1}_{\{\|G^{-1}(x-m)\|\leq r\}
\cap
B_R}\dvt
m\in\R^k, G\in\operatorname{PDS}(k), r>0\bigr\},
\end{eqnarray*}
have uniformly bounded envelopes. According to Lemma 3 in
\cite{lopuhaa97}, the corresponding classes of graphs have polynomial
discrimination. Therefore, Theorem \ref{th:pollard} also applies to the
components of $\Psi_1$ and $\Psi_2$. It follows that
\[
\label{eq:expansion theta} 0 = \Lambda(\widehat{\theta}_n)+\int
\Psi(y,\theta_0) (P_n-P_0)(\mathrm{d}y)+\mathrm{o}_\mathbb{P}(n^{-1/2}).
\]
Now, $\Lambda(\theta_0)=0$ and since $\Psi(y,\theta_0)$ has bounded
support, the term
\[
\int\Psi(y,\theta_0) (P_n-P_0)(\mathrm{d}y) = \frac1n\sum_{i=1}^n
\bigl(\Psi(X_i,\theta_0)-\mathbb{E}\Psi(X_i,\theta_0)\bigr),
\]
behaves according to the central limit theorem and is therefore of the
order $\mathrm{O}_\mathbb{P}(n^{-1/2})$. Because $\widehat{\theta}_n\to
\theta_0$
with probability one, according to Theorem \ref{th:consistency
estimator}, we find
\[
0=\Lambda'(\theta_0)(\widehat{\theta}_n-\theta_0)+\mathrm{O}_\mathbb
{P}(n^{-1/2})+\mathrm{o}_\mathbb{P}(\|\widehat{\theta}_n-\theta_0\|).\vadjust{\goodbreak}
\]
Because $\Lambda'(\theta_0)$ is non-singular, this gives
$\|\widehat{\theta}_n-\theta_0\|=\mathrm{O}_\mathbb{P}(n^{-1/2})$ and when
inserting this, we conclude that
\[
\Lambda'(\theta_0)(\widehat{\theta}_n-\theta_0)=-\frac1{n}\sum_{i=1}^n
\bigl(\Psi(X_i,\theta_0)-\mathbb{E}\Psi(X_i,\theta_0)\bigr) +
\mathrm{o}_\mathbb{P}(n^{-1/2}),
\]
which proves the first statement. For the second statement note that
\[
\int\bigl(\mathbh{1}_{E_n}(y)-\phi_n(y)\bigr)\|y\|^2 P_n(y)=\mathrm{O}(n^{-1}),
\]
with probability one. This follows from the characterization given in
Theorem \ref{th:characterization} and the fact that $P_0$ satisfies
\eqref{eq:atoms}. This means that the MCD functional $\theta_n$ also
satisfies equation~\eqref{eq:expansion emp}. From here on the argument
is the same as before, which proves the theorem.
\end{pf*}

\begin{pf*}{Proof of Theorem \ref{th:IF}}
Consider expansion
\eqref{eq:IF expansion theta} and write $E_0=E(\mu_0,\Sigma_0,\rho_0)$.
Because, according to Theorem \ref{th:continuity},
$(\mu_{\eps,x},\Gamma_{\eps,x},\rho_{\eps,x})\to
(\mu_0,\Gamma_0,\rho_0)$, as $\eps\downarrow0$, for $x\notin
\partial
E_0$ we get $\phi_{\eps,x}(x)\to\mathbh{1}_{E_0}(x)$ and hence
\[
\lim_{\eps\downarrow0}\Phi_{\eps}(x)=\Psi(x,\theta_0),
\]
with $\Psi$ defined in \eqref{eq:def Psi}. Because
$\mathbh{1}_{E_0^\circ}\leq\phi_0 \leq\mathbh{1}_{E_0}$, it follows
from \eqref{eq:mass En} that $\Lambda(\theta_0)=0$. As~$\Lambda$
has a
non-singular derivative at $\theta_0$, we find from \eqref{eq:IF
expansion theta},
\[
0 = (1-\eps) \Lambda'(\theta_0)(\theta_{\eps,x}-\theta_0) + \eps
\Phi_{\eps}(x) +\mathrm{o}(\|\theta_{\eps,x}-\theta_0\|),
\]
from which we first deduce that $\theta_{\eps,x}-\theta_0=O(\eps)$, and
then obtain the influence function:
\[
\operatorname{IF}(x;\Theta,P_0) =
\lim_{\eps\downarrow0}\frac{\theta_{\eps,x}-\theta_0}{\eps} =
-\Lambda'(\theta_0)^{-1}\Psi(x).
\]
\upqed\end{pf*}
\end{appendix}

\section*{Acknowledgements}
We like to thank two anonymous referees for their comments and suggestions,
which were helpful for improving the manuscript.
Special thanks goes to referee~1 for correcting a~mistake in the proof
of Corollary 4.1
and for pointing out the relation with other relevant
work \cite
{gordaliza1991,cuesta-albertosgordalizamatran1995,cuesta-albertosgordalizamatranAS1997,garcia-escuderogordalizamatranAS1999},
which has lead to shorter proofs of Proposition \ref{prop:smalleig}
and the first part of Theorem \ref{th:characterization}.

% imsref loaded by svajune.rapalyte, 2011-10-05 10:42:11
%

\printhistory

\end{document}